\documentclass{amsart}
\usepackage[british]{babel}
\usepackage[OT2,T2A,T1]{fontenc}
\usepackage[utf8]{inputenc}		
\usepackage{eulervm}	
\usepackage[scr=dutchcal, bb=boondox]{mathalfa}

\usepackage{float}
\usepackage{extarrows}

\usepackage{titlesec}
\titleformat{\section}
	{\centering\Large\scshape\bfseries}{\thesection.}{1em}{} 
\titleformat{\subsection}
 	{\normalfont\bfseries\scshape}{\thesubsection.}{.5em}{}
\titleformat*{\paragraph}{\bfseries}

\titleformat{\subsubsection}
 	{\normalfont\bfseries\scshape}{\thesubsubsection.}{.5em}{}
\titleformat*{\paragraph}{\bfseries}

\usepackage{titletoc}

    \titlecontents{section}
    [0em] %
    {\medskip}
    {\thecontentslabel\quad}
    {}
    {\hfill\contentspage}

\usepackage{amsmath,amsthm,amssymb,dsfont,extarrows,amscd,amsfonts,mathtools,bm,physics}
\usepackage[v1]{subfiles}
\usepackage{hyperref}
\usepackage{mathrsfs}
\usepackage[mathscr]{euscript}
\usepackage{tensor}
\usepackage[most]{tcolorbox}
\usepackage{stmaryrd}
\usepackage[capitalize]{cleveref}
\usepackage{MnSymbol}
\usepackage{thmtools}

\usepackage{tikz}
\usepackage{tikz-cd}
\usetikzlibrary{graphs,decorations.pathmorphing,decorations.markings}
\tikzset{
    labl/.style={anchor=south, rotate=90, inner sep=.5mm}
}

\usepackage{todonotes}
\usepackage{geometry}
\usepackage[alphabetic, initials]{amsrefs}

\usepackage{titlesec}

\titleformat{\paragraph}
{\normalfont\normalsize\bfseries}{\theparagraph}{1em}{}
\titlespacing*{\paragraph}
{0pt}{2.5ex}{1ex}

\DeclareFontFamily{U}{MnSymbolC}{}
\DeclareSymbolFont{MnSyC}{U}{MnSymbolC}{m}{n}
\DeclareFontShape{U}{MnSymbolC}{m}{n}{
	<-6>  MnSymbolC5
	<6-7>  MnSymbolC6
	<7-8>  MnSymbolC7
	<8-9>  MnSymbolC8
	<9-10> MnSymbolC9
	<10-12> MnSymbolC10
	<12->   MnSymbolC12}{}
\DeclareMathSymbol{\intprod}{\mathbin}{MnSyC}{'270}

\usepackage{stmaryrd}
\SetSymbolFont{stmry}{bold}{U}{stmry}{m}{n}

\theoremstyle{plain} 
\newtheorem{thm}{Theorem} [section]
\newtheorem{cor}{Corollary}[section]
\newtheorem{lem}{Lemma} [section]
\newtheorem{prop}{Proposition}[section]

\theoremstyle{definition}
\newtheorem{defn}{Definition}[section]
\newtheorem{ex}{Example}[section]
\newtheorem{conj}{Conjecture}[section]

\newtheorem{oss}{Remark}[section]

\numberwithin{equation}{section}

\newcommand{\CC}{\mathbb{C}}

\DeclareMathOperator{\X}{\mathfrak{X}}
\DeclareMathOperator{\Lie}{\mathscr{L}}

\newcommand{\Hol}{\mathscr{O}}

\newcommand{\EE}{\mathbb{E}}

\newcommand{\FM}{Frobenius manifold~}
\newcommand{\gFM}{generalized Frobenius manifold~}

\newcommand{\gFMs}{generalized Frobenius manifolds~}

\title[Generalized Frobenius Manifolds and their symmetries]{Generalized Frobenius Manifolds: symmetries and their classification in 2-dimensions}

\author[A.~Bedell]{A.~Bedell}
\author[M.~Feigin]{M.~Feigin}
\author[L.~Kaminski]{L.~Kaminski}
\author[A.~Proserpio]{A.~Proserpio}
\author[I.~A.~B.~Strachan]{I.~A.~B.~Strachan}

\address{School of Mathematics and Statistics,
University of Glasgow, Glasgow, G12 8QQ, United Kingdom.}

\email{andrebedell@gmail.com}
\email{misha.feigin@glasgow.ac.uk}
\email{ale.prose.maths@gmail.com}
\email{ljmkaminski@gmail.com}
\email{ian.strachan@glasgow.ac.uk}

\begin{document}

\renewcommand{\hbar}{\hslash}

\begin{abstract}
A theory of symmetries of generalized Frobenius manifolds is developed. This is used to study the solution space of all generalized Frobenius manifolds that are \lq above\rq~a specific Frobenius manifold. Since all (semi-simple) Frobenius manifolds in 2-dimensions are known, by constructing all homogeneous Legendre fields associated to these manifolds we classify all 2-dimensional generalized Frobenius manifolds. With this description, the known examples of generalized Frobenius manifolds correspond to distinguished points in this solution space. For Hurwitz-Frobenius manifolds the action of a Legendre transformation on Hurwitz space data, and in particular the primary differential, is given.
\end{abstract}

\maketitle

\tableofcontents


\section{Introduction}

The Witten-Dijkgraaf-Verlinde-Verlinde (or WDVV) equations of associativity (here, and in what follows, all indices run over $1\,,\ldots\,,n$)
\begin{equation}
\frac{\partial^3F(t)}{\partial t^\alpha\partial t^\beta\partial t^\nu}\eta^{\nu\mu}
\frac{\partial^3F(t)}{\partial t^\mu\partial t^\gamma\partial t^\lambda}-
\frac{\partial^3F(t)}{\partial t^\beta\partial t^\gamma\partial t^\nu}\eta^{\nu\mu}
\frac{\partial^3F(t)}{\partial t^\alpha\partial t^\mu\partial t^\lambda}=0
\label{wdvv}
\end{equation}
were originally introduced in the study of topological quantum field theory \cites{wdvv1,wdvv2}. Geometrically, a solution defines a commutative, associative multiplication $\circ\,:  \X_M \times \X_M\rightarrow  \X_M$ of vector fields, i.e.
\[
\frac{\partial~}{\partial t^\alpha} \circ \frac{\partial~}{\partial t^\beta} = \eta^{\mu\nu}\frac{\partial^3 F}{\partial t^\alpha \partial t^\beta \partial t^\mu} \,\frac{\partial~}{\partial t^\nu}\
\]
satisfying the Frobenius property $\eta(X\circ Y,Z)=\eta(X,Y\circ Z)\,.$
The metric $\eta$, used to raise and lower indices, is flat and the coordinates $\{ t^\alpha \}$ are flat coordinates, i.e. the components of the metric is this coordinate system are constants. The WDVV equations are the conditions for this multiplication to be associative.

The WDVV equations can be supplemented with additional conditions and structures. A (Dubrovin-)Frobenius manifold is one such example of this - here one has a covariantly constant unity field as well as an Euler field which generates scaling symmetries \cite{Dubrovin1996}. Applications of such a structure abound - Frobenius manifolds play a central role in the Gromov-Witten theory, in singularity theory and the bi-Hamiltonian structure of certain integrable hierarchies.

But it was soon found that this was too rigid a structure - other theories had the WDVV equations as a central component but were not Frobenius manifolds. For example, in quantum K-theory one has the WDVV equations with the unity field which is not covariantly constant \cites{K1,K2}. Further examples followed. In the theory of Hurwitz spaces one requires a so-called primary differential in order to define the multiplication and metric and in \cite{Dubrovin1996} those (termed admissible differentials) which produced a Frobenius manifold were classified. But there are others. Motivated by the bi-Hamiltonian structure of the Ablowitz-Ladik (dispersionless) hierarchy, the solution to the WDVV equations 
\begin{equation}
F=\frac{1}{2} (t^1)^2 t^2 + t^1 e^{t^2} + \frac{1}{2} (t^1)^2 \log t^1
\label{nonFMsolution}
\end{equation}
with Euler field
\[
E=t^1 \frac{\partial~~}{\partial t^1} + \frac{\partial~~}{\partial t^2}\,
\]
and a non-flat unity field
\[
e=\frac{1}{t^1-e^{t^2}} \left( t^1 \frac{\partial~~}{\partial t^1} + \frac{\partial~~}{\partial t^2}\right)
\]
was obtained in \cite{BCR}, together with a Hurwitz space description using a non-admissible differential. This idea was further developed in \cite{RomanoThesis} and \cite{BCRR} and such results were one of the motivations behind the formal definition of a generalized Frobenius manifold \cite{GF1} (see Section \ref{gFMdef}).

The same structures can be obtained via another route by studying symmetries of the WDVV equations. A symmetry of the WDVV equations are transformations
\begin{itemize}
\item[(i)] $t^\alpha \mapsto {\tilde t}^\alpha\,,$
\item[(ii)] $\eta_{\alpha\beta} \mapsto {\tilde\eta}_{\alpha\beta}\,,$
\item[(iii)] $F \mapsto {\tilde F}$
\end{itemize}
that preserve the WDVV equations, i.e. which map solutions of the WDVV equations to new solutions. In \cite{Dubrovin1996} two classes of symmetries were constructed:

\begin{def}\label{symdef}
{}
\begin{itemize}

\item[(a)]\emph{Legendre-type transformations $S_\kappa$}: This is defined by
\begin{eqnarray*}
{\hat t}_\alpha & = & \partial_{t^\alpha} \partial_{t^\kappa} F(t)\,,
\qquad ({\rm note:}\quad{\hat t}_\alpha=\eta_{\alpha\beta}{\hat t}^\beta)\\
\frac{\partial^2{\hat F}}{\partial{\hat t}^\alpha \partial{\hat t}^\beta} & = &
\frac{\partial^2{     F}}{\partial{     t}^\alpha \partial{     t}^\beta} \,,\\
{\hat\eta}_{\alpha\beta} & = & \eta_{\alpha\beta}\,.
\end{eqnarray*}

\noindent Here $\kappa$ is a distinguished index which labels the Legendre transformation. We assume that the vector field $\partial_{t^\kappa}$ is invertible.

\medskip

\item[(b)]\emph{Inversion Symmetry $I$:} This is defined by
\begin{eqnarray*}
\hat{t}^1  & = & \frac{1}{2}\frac{t_{\sigma}t^{\sigma}}{t^n}, \\ \hat{t}^{\alpha}  & = & \frac{t^{\alpha}}{t^n}, \quad\mbox{  (for }\alpha \neq 1, n), \\ \hat{t}^n  & = & -\frac{1}{t^n},
\end{eqnarray*}
and the transformation
\begin{eqnarray*}
\hat{F}(\hat{t}) & = & (\hat{t}^n)^2 F \left( \frac{1}{2}\frac{\hat{t}_{\sigma}\hat{t}^{\sigma}}{\hat{t}^n}, -\frac{\hat{t}^{2}}{\hat{t}^n}, ... , -\frac{\hat{t}^{n-1}}{\hat{t}^n},
-\frac{1}{\hat{t}^n} \right) + \frac{1}{2}\hat{t}^1\hat{t}_{\sigma}\hat{t}^{\sigma},\\
\hat{\eta}_{\alpha\beta} & = & \eta_{\alpha\beta}.
\end{eqnarray*}
\end{itemize}
\end{def}

\medskip

\noindent Moreover, these transformations are more than just a map between solutions of the WDVV equations: they preserves the stronger notion of a Frobenius manifold. Furthermore, the monodromy data is preserved under such a transformation, as well as (in the semi-simple case) Stokes and other data \cite{yang}.


The Legendre transformation may be reinterpreted more geometrically in terms of a  vector field \cites{RS,SS}
\[
\delta = \frac{\partial~}{\partial t^\kappa}\,, \qquad\qquad \kappa\in\{1\,,\ldots\,,n \}\,.
\]
With this the metric transformation may be written as
\[
{\tilde\eta}(X,Y) = \eta(\delta\circ X,\delta\circ Y)
\]
for arbitrary vector fields $X$ and $Y\,.$ The field $\delta$ is distinguished by the fact that it is flat for the Levi-Civita connection of the metric $\eta\,,$ namely $\nabla_X\delta=0\,$ and homogeneous with respect to the Euler field of the manifold.

\begin{figure}[ht]
\begin{tikzpicture}
\draw[thick] (0,0) parabola (1.5,2); 
\draw[thick,->] (0.01,0) -- (0,0);
\node[] at (-0.5,0.02) {FM $\cross$};
\node[] at (1.95,2.2) {$\cross$ gFM};
\node[] at (3.2,1) {$\delta$-flat and homogeneous};
\end{tikzpicture}
\caption{A flat, homogeneous, Legendre field maps a generalized Frobenius manifold to a Frobenius manifold.}
\label{fig:gfm2fm}
\end{figure}
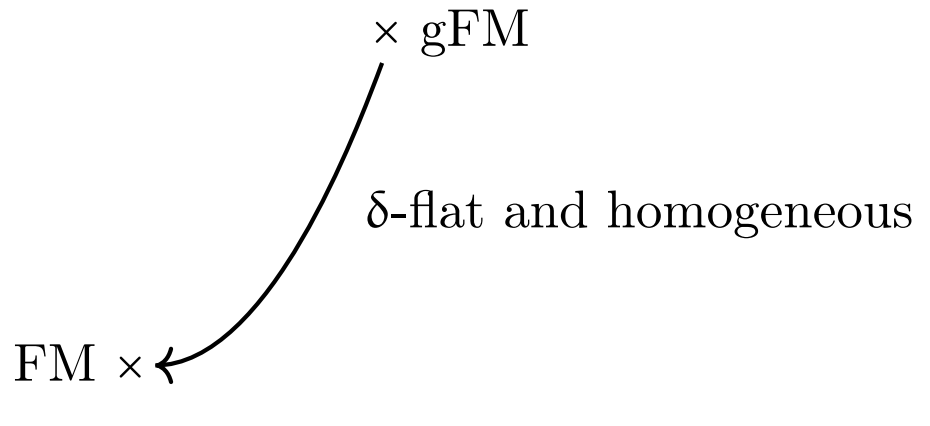

In \cite{antidiag} it was shown that applying a Legendre transformation to the solution (\ref{nonFMsolution}) results in the solution
\[
F=\frac{1}{2} (t^1)^2 t^2 + e^{t^2}
\]
which {\sl does} define a Frobenius manifold. This holds more generally - and in \cite{LegendreGFM} it was observed that applying such a Legendre transformation to a generalized Frobenius manifold maps it to a Frobenius manifold rather than another generalized Frobenius manifold. This is shown schematically in Figure \ref{fig:gfm2fm}.

\noindent One may construct Legendre transformations with more general vector fields than the coordinate vector fields $S_\kappa\,.$ In \cite{SS} the notion of a Legendre field was developed. Such fields\footnote{These were called generalized Legendre fields in \cite{SS}. We now just refer to these as Legendre fields, with fields that are covariantly constant as flat Legendre fields.} satisfy the defining relation
\begin{equation}
X\circ \nabla_Y\delta = Y\circ \nabla_X\delta
\label{def1}
\end{equation}
and map solutions of the WDVV equations to new solutions, i.e. they generate a symmetry of the WDVV equations. It was also proved in \cite{SS} that the inverse field $\delta^{-1}$ (assuming it exists) is also a Legendre field and that there is a transitive action under composition. As remarked in \cite{SS} non-flat Legendre fields do not preserve the condition $\nabla e=0$ and hence do not map Frobenius manifolds to Frobenius manifolds. They are now best thought of as defining symmetries of generalized Frobenius manifolds (see Section \ref{LFsection}). 

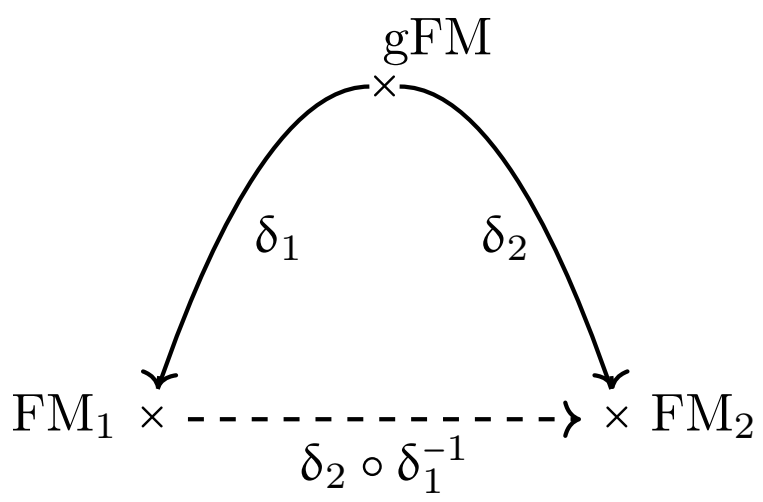
\begin{figure}[ht]
\begin{tikzpicture}
\draw[thick] (1.4,2) parabola (0,0) ; 
\draw[thick,->] (0.01,0.085) -- (0,0);
\draw[thick] (1.6,2) parabola (3,0) ; 
\draw[thick,->] (3,0.085) -- (3.01,0);
\node[] at (-0.45,-0.18) {FM$_1$ $\cross$};
\node[] at (3.45,-0.18) {$\cross$ FM$_2$};
\node[] at (1.5,2.0) {$\cross$};
\node[] at (1.85,2.3) {gFM};
\node[] at (0.8,1) {$\delta_1$};
\node[] at (2.3,1) {$\delta_2$};
\draw[thick, dashed,->] (0.2,-0.2) -- (2.8,-0.2);
\node[] at (1.5,-0.5) {$\delta_2 \circ \delta_1^{-1}$};
\end{tikzpicture}
\caption{Transitive action: different flat, homogeneous Legendre fields from a generalized Frobenius manifold generate different Frobenius manifolds which are connected by a flat, homogeneous Legendre field.}
\label{fig:transitive}
\end{figure}

Putting these results together give the following: starting from a specific generalized Frobenius manifold, there are many possible flat, homogeneous, Legendre fields, each of which will generate a Frobenius manifolds. Using the transitive action (see Proposition \ref{transitive}) under composition, these Frobenius manifolds are themselves connected by a flat, homogeneous, Legendre field (see Figure \ref{fig:transitive}). Inverting this construction shows that two different Frobenius manifolds may be Legendre-transformed to the {\sl same} generalized Frobenius manifold.

\noindent Alternatively, any specific generalized Frobenius manifold will be \lq above\rq~different Frobenius manifolds that are themselves related by flat, homogeneous Legendre transformations (for a specific example of this construction, see Example \ref{ex3.2}).

\begin{figure}[ht]
\begin{tikzpicture}
\draw[->] (-0.2,0.1) -- (-0.2,3.5);
\draw (-0.2,1.5) -- (1.0,2.7) -- (4.5,2.7) -- (3.3,1.5) -- cycle;
\draw[thick,->] (0,0) parabola (1.5,2); 
\node[] at (-0.53,0.02) {FM $\cross$};
\node[] at (1.95,2.2) {$\cross$ gFM};
\node[] at (1.4,1) {$\delta$};
\node[] at (-0.5,3.5) {${\tilde{d}}$};
\node[] at (5,2.8) {$\mathbb{P}^{n-1}_{d,\tilde{d}}$};
\end{tikzpicture}
\caption{The space of generalized Frobenius manifolds that can be constructed from a given Frobenius manifold of charge $d$ with a homogeneous Legendre field.}
\label{fig:fm2gfm}
\end{figure}
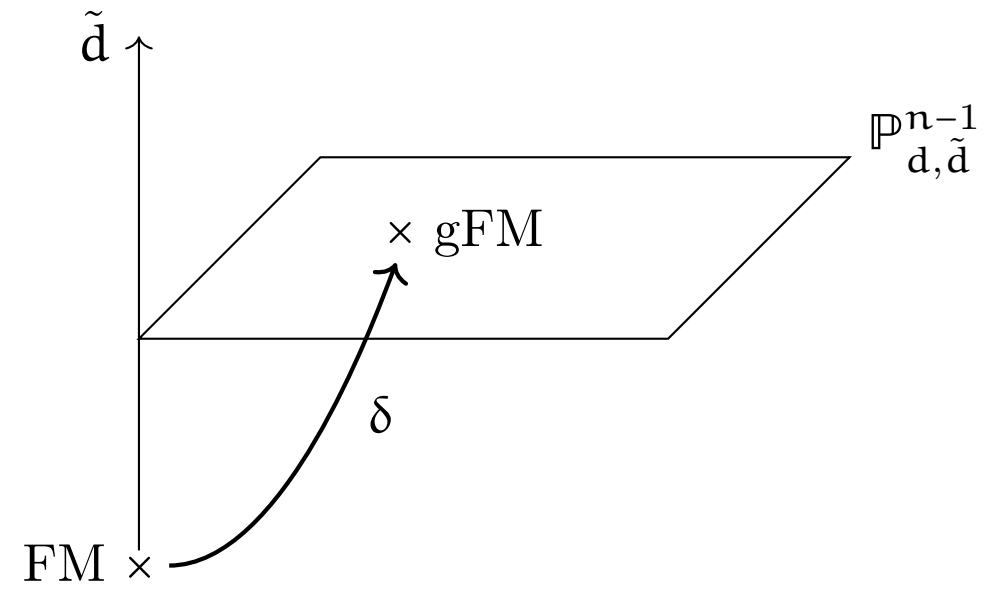

The purpose of this paper is two-fold. Firstly, to develop a full theory of symmetries of generalized Frobenius manifolds, including inversion symmetries, and secondly to start with a \lq base\rq~Frobenius manifold and construct Legendre transformations to the space of generalized Frobenius manifolds. This will, in general (and as observed in \cite{SS}), not preserve the flatness of the unity field $e$ and hence will map to a generalized Frobenius manifold rather than to a Frobenius manifold. But given a Frobenius manifold, the space of homogeneous Legendre fields depends on data: one has the charge of the Frobenius manifold, the homogeneity of the Legendre field (homogeneity being required to preserve the properties of the Euler field) and as defining equations (\ref{def1}) for a Legendre field are linear differential equation, solutions will depend on constants of integration. This produces, for a given Frobenius manifold of charge $d$, an $(n-1)$-dimensional complex space $\mathbb{P}_{d,\tilde{d}}^{n-1}$ (isomorphic to the complex projective space $\mathbb{CP}^{n-1}$) of Legendre fields that generate \gFMs of charge $\tilde{d}$ which have the same underlying Frobenius manifold. The new charge $\tilde{d}$ depends on the original charge $d$ and the homogeneity of the Legendre field. This is shown schematically in Figure \ref{fig:fm2gfm}.

Since two-dimensional semi-simple Frobenius manifold have been classified, one can, by constructing all their homogeneous Legendre fields, classify all 2-dimensional generalized Frobenius manifolds. To reduce this problem, we introduce the definition of a canonical \gFM. This fixes the charge of the \gFM to be $\tilde{d}=1\,.$ Since Legendre fields can be composed, with a transitive action, any \gFM may be transformed by a suitable Legendre transformation to a canonical generalized Frobenius manifold. There are other geometric reasons why canonical \gFMs are distinguished (see Section \ref{classification}), and we also note that all known examples of \gFMs have $\tilde{d}=1\,.$

\begin{figure}[ht]
\begin{tikzpicture}
\draw[->] (-0.2,0.1) -- (-0.2,3.5);
\draw (-0.2,1.5) -- (1.0,2.7) -- (5.5,2.7) -- (4.3,1.5) -- cycle;
\draw[thick,->] (0,0) parabola (1.4,2); 
\draw[thick,->] (0,0) parabola (3.4,2); 
\node[] at (-0.70,0.02) {$\mathbb{C}^{2}/A_2\cross$};
\node[] at (1.96,2.2) {$\cross\,\,\mathbb{C}^2/{A_2^{aff}}$};
\node[] at (3.96,2.2) {$\cross\,\,\mathbb{C}^2/{G_2^{aff}}$};
\node[] at (1.25,1) {$\delta_1$};
\node[] at (2.8,1) {$\delta_2$};
\node[] at (-0.5,3.5) {${\tilde{d}}$};
\node[] at (6.3,2.8) {$\mathbb{P}^{1}_{d=\frac{1}{3},\tilde{d}=1}$};
\node[] at (-0.5,1.5) {${}_{\tilde{d}=1}$};
\end{tikzpicture}
\caption{Two \gFMs can be related to the same underlying Frobenius manifold.}
\label{fig:fig4}
\end{figure}
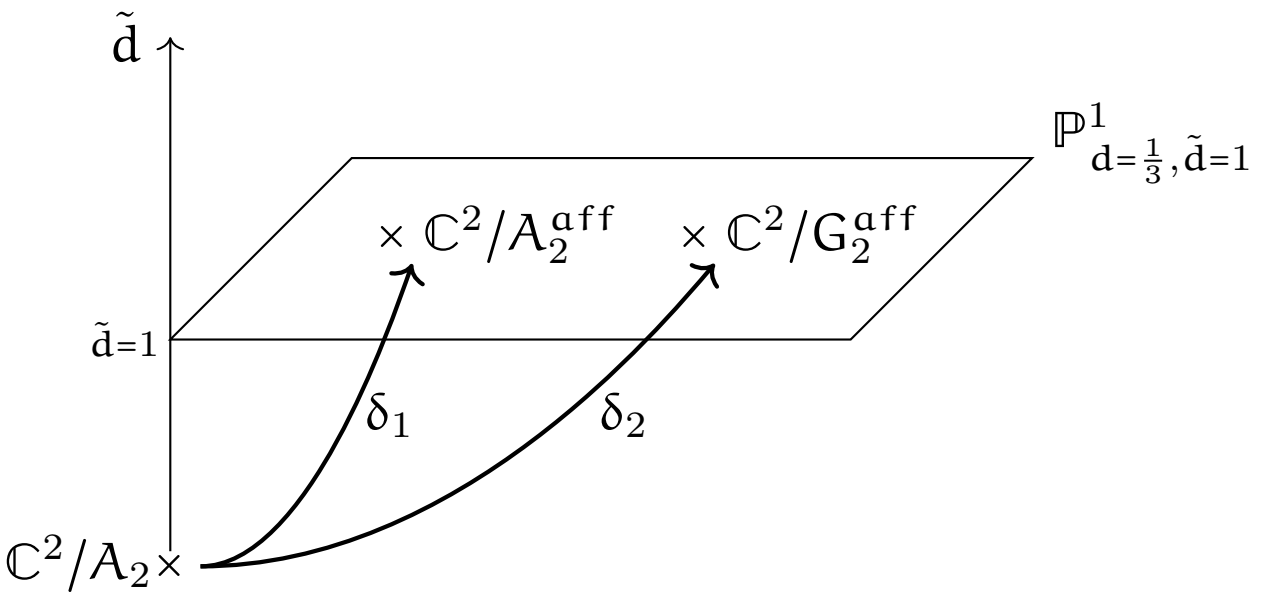

Even after fixing the homogeneity of the Legendre field to give a canonical generalized Frobenius manifold, one has a space of Legendre fields, and in examples, certain points are distinguished by the algebraic properties of the solution. As an example of this, the two of the \gFMs constructed in \cite{orbit1} via an affine orbit space construction have the same charge and the same Euler field, and can both be mapped to the same underlying Frobenius manifold, namely 
the Frobenius manifold on the orbit space
$\mathbb{C}^2/A_2\,.$ Inverting this, one has two Legendre 
fields (denoted $\delta_1$ and $\delta_2$ in Figure \ref{fig:fig4}) of the same homogeneity and satisfying the same defining equation. What distinguishes these fields from each other, and from others in the space of all such fields, are certain algebraic properties. Since Legendre fields are assumed to have inverses and can be composed, one can also construct a symmetry, generated by the Legendre field $\delta_1^{-1}\circ\delta_2$ (or its inverse), between these two generalized Frobenius manifolds.

The rest of the paper is laid out as follows: Section 2 gives the definition of a generalized Frobenius manifold and other key ideas before Legendre transformations and Legendre fields are defined in Section 3. In this section  the theory is further developed and a canonical class of Legendre fields is introduced. In Section 4, building on the classification of $2D$ Frobenius manifolds (see Table 1), all homogeneous Legendre fields are explicitly constructed in terms of solutions of the hypergeometric equation, and this provides explicit maps from a specific Frobenius manifold to a space of generalized Frobenius manifolds. These results are contained in Tables 1 and 2. In Section 5 this is reversed: the known $2D$ examples of generalized Frobenius manifolds are placed within this wider construction. Section 6 starts with a study of associativity equations in $2D$ and the solution space is studied, including a phase space analysis. Section 7 explores the interaction between inversion symmetries and Legendre symmetries, and this is then used to construct an inversion symmetry for generalized Frobenius manifolds. An important class of Frobenius and generalized Frobenius manifolds exists on Hurwitz spaces - moduli spaces of holomorphic maps. Central to these constructions are a pair $(\lambda,\omega)$ where $\lambda$ is the Landau-Ginzburg superpotential and $\omega$ the primary differential, or Saito form. The main result in this section is the action of a Legendre transform on this data. The paper concludes with some further discussion and potential applications.

\section{Preliminaries}\label{gFMdef}

Following {\cite{GF1}}, we define: 
\begin{defn}[Generalized Frobenius manifold]
An $n$-dimensional generalized Frobenius manifold of charge $d\in\CC$ is an $n$-dimensional complex manifold\footnote{The definition also works in the category of smooth manifolds. One just has to replace “holomorphic” with “smooth”.} $M$ equipped with:
\begin{itemize}
    \item an $\Hol_M$-bilinear multiplication $\circ:\X_M\otimes_{\Hol_M}\X_M\to \X_M$ on the sheaf of holomorphic tangent vectors, called \emph{Frobenius product}.
    \item a symmetric and non-degenerate bilinear form $\eta:\X_M\otimes_{\Hol_M}\X_M\to\Hol_M$, called \emph{Frobenius pairing} (or metric).
    \item two holomorphic vector fields $e$ and $E$, respectively called \emph{unity} and \emph{Euler vector field}.
\end{itemize}
such that:\begin{enumerate}
    \item [(i)] At any $p\in M$, $T_pM$ equipped with the multiplication $\circ_p$ and the non-degenerate, $\CC$-bilinear form $\eta_p$ is a Frobenius algebra.
    \item [(ii)] $e_p$ is the identity of $T_pM$ at any $p\in M$.
        \item [(iii)] for any $\hbar \in\CC$, the \emph{Dubrovin connection} ${}^\hbar\nabla:={}^\eta\nabla+\hbar\,\circ$ is flat and torsion-free, ${}^\eta\nabla$ being the Levi-Civita connection of $\eta$.
    \item [(iv)] $E$ satisfies the following relations:
      \begin{eqnarray*}
\Lie_E\circ& = &\circ \, ,\\
\Lie_E\eta\,& = & \, (2-d)\,\eta\,,
    \end{eqnarray*}
    
\noindent where $\Lie$ is the Lie derivative.
\end{enumerate}
\end{defn}

\begin{oss} The vanishing of the torsion of the Dubrovin connection is equivalent to the commutativity of the product. The curvature has three terms: the $O(\hbar^0)$ implies the curvature of the connection ${}^\eta\nabla$ is zero and so the metric is flat and hence there exists a set of flat coordinates $\{ t^i\,,i=1\,,\ldots\,,n \}$ in which the components of the metric are constants; the $O(\hbar^1)$ implies that there exists a scalar prepotential $F$ and the $O(\hbar^2)$ terms gives associativity and hence the WDVV equations for the prepotential (\ref{wdvv}).

\end{oss}

From condition $(iv)$ the Euler field is a conformal/homothetic Killing vectors and since the metric is flat this must be linear in the flat coordinates (excluding the inversions in infinity which was mentioned in \cite{Dubrovin1996} but which has not been developed further). Thus, over redefinitions of coordinates, this can be taken to be
\[
E=\sum_{i=1}^n d^i t^i \frac{\partial~~}{\partial t^i} + \sum_{ \{ i|d^i=0\} } r^i \frac{\partial~~}{\partial t^i}
\]
where the $d^i$ and the $r^i$ are constants. The spectrum 
$\{\mu^i\}$ of the generalized Frobenius manifold is expressed in terms of the endomorphism
\begin{equation}
\mu(X) = \frac{2-d}{2} \, X - \nabla_XE\,,
\label{endo}
\end{equation}
evaluated in the flat coordinates of the metric, and hence
\begin{equation}
d^i = 1 - \left(\mu_i + \frac{d}{2}\right)\,.
\label{spectrum}
\end{equation}
A standard calculation gives the endomorphism $\mu$ is skew-adjoint with respect to the metric,
\begin{equation}
\eta(\mu(X),Y) + \eta(X,\mu(Y)) = 0\,
\label{skew}
\end{equation}
and (assuming the metric is antidiagonal, with $\eta_{ij} = \delta_{i+j,n+1}$) this implies the
condition
\begin{equation}\label{muconstraint}
\mu_i + \mu_{n+1-i} = 0 \,.
\end{equation}
All this has been developed without the need of the unity field to be covariantly constant.

\begin{defn}
    A Frobenius manifold is a generalized Frobenius manifold with the additional condition $\nabla_Xe=0\,$ for all vector fields $X\,.$
\end{defn}

\noindent Semi-simple Frobenius manifolds in 2-dimensions have been classified, and these are given (for the cases where $\eta(e,e)=0$) in Table 1. If $\eta(e,e)\neq 0$ then the prepotential is cubic.

\section{Legendre fields}\label{LFsection}

Symmetries of the WDVV equations of Legendre type can be generated by more general vector fields than just coordinate vector fields. This was introduced in \cite{SS} and developed further in \cites{GF1,LegendreGFM, LeoThesis}. Here we draw together these results and extend them.

\subsection{Legendre fields and their properties}

Consider the following structure: $ \{ M, \eta,\nabla,\circ\}$, where $M$ is a manifold and
\begin{itemize}
\item[(a)] $\eta$ is a metric;
\item[(b)] $\nabla$ is a connection (not, at this stage, metric or torsion free);
\item[(c)] $\circ: \X_M \times \X_M \rightarrow \X_M$ is a commutative, associative multiplication with a unit $e$ which is compatible with the metric, i.e.
\[
\eta(X \circ Y,Z) = \eta(X, Y\circ Z)\,;
\]
\item[(d)] $\nabla\circ$ is a completely symmetric 
$(1,3)$-tensor.
\end{itemize}
A Legendre field will map metrics and connections to new metric and connections which preserve certain basic properties.
\begin{defn}
Let $\delta$ be an invertible vector field (so $\delta^{-1} \circ \delta = e$). Let
\begin{eqnarray*}
{\tilde\eta}(X,Y) & = & \eta( \delta\circ X,\delta\circ Y)\,,\\
{\tilde\nabla}_XY&=& \delta^{-1} \circ \nabla_X(\delta\circ Y)\,.
\end{eqnarray*}
\end{defn}
\noindent The following follows from direct calculation.
\begin{lem}\label{firstlemma}
For all vector fields $X\,,Y\,,Z \in \X_M$
\begin{eqnarray*}
( {\tilde\nabla}_X{\tilde{\eta}})(Y,Z) & = & \left( \nabla_X\eta\right)(\delta\circ Y,\delta\circ Z)\,,\\
T^{\tilde\nabla}(X,Y) & = & T^\nabla(X,Y) + \delta^{-1} \circ \left\{ Y \circ \nabla_X \delta - X \circ \nabla_Y \delta\right\}\,,\\
R^{\tilde\nabla}(X,Y)(Z) & = & \delta^{-1} \circ R^\nabla(X,Y) (\delta\circ Z)\,,
\end{eqnarray*}
where $T^\nabla$ and $R^\nabla$ are the torsion and curvature tensors of the connection $\nabla\,.$
\end{lem}
\noindent This result motivates the following \cite{SS}:
\begin{defn}
A (generically) invertible vector field $\delta \in \X_M$ is a Legendre field if it satisfies the equation
\begin{equation}\label{defLF}
X \circ \nabla_Y \delta = Y \circ \nabla_X \delta
\end{equation}
for all $X\,,Y\in \X_M\,.$
\end{defn}
The definition appears to be overdetermined, with more equations than unknowns. However if one sets $Y=e$ one obtains the manifestly determined system
\[
\nabla_X\delta = X \circ \nabla_e\delta\,.
\]
Conversely, a solution of this equation implies that $\delta$ is a Legendre field, so this provides an alternative definition of a Legendre field. Thus a Legendre field will map torsion free connections to torsion free connections, metric connections to metric connections, and zero-curvature conditions to zero-curvature conditions. Since the WDVV equations are equivalent to the vanishing of the curvature of Dubrovin connection, by Lemma \ref{firstlemma} a Legendre field will generate a new  - Legendre transformed - solution of the WDVV equations. It is important to note that geometrically the multiplication does not change when such a transformation is applied. However, the coordinate expression for this tensor does depend on the choice of flat-coordinate system and these do change under such a transformation. We now assume that $\{M\,,\eta\,,\nabla\,,\circ\}$ together with identity vector field $e$ and Euler vector field $E$ defines generalized Frobenius manifold. 

Legendre fields are always defined with respect to a specific connection. So while $\delta^{-1}$ can be defined algebraically, it turns out geometrically that this is a Legendre field with respect to the connection $\widetilde{\nabla}$, not the connection $\nabla.$ Similarly, Legendre fields can be combined if one takes care to track which connections are being used. This is encapsulated in the following \cite{SS}.

\begin{prop}\label{transitive}

\phantom{text}

\begin{itemize}
\item[(i)] The unity field $e$ is a Legendre field.
\vskip 2mm
\item[(ii)] If $\delta$ is a Legendre field for $\nabla$ then $\delta^{-1}$ is a Legendre field for $\tilde\nabla\,.$
\vskip 2mm
\item[(iii)] Consider connections $\nabla^{(i)}\,,i=1,2,3$ with interconnecting Legendre fields $\delta^i_j\,,$
\[
\nabla^{(1)} \xrightarrow{~~~\delta^2_1~~~} \nabla^{(2)}\xrightarrow{~~~\delta^3_2~~~} \nabla^{(3)}\,.
\]
Then the product $\delta^3_1 = \delta^2_1 \circ \delta^3_2$ is a Legendre field for the connection $\nabla^{(1)}$ directly connecting to the connection $\nabla^{(3)}\,:$
\[
\nabla^{(1)} \xrightarrow{~~~\delta^3_1~~~} \nabla^{(3)}
\]

\end{itemize}
\end{prop}

\noindent This gives a transitive, groupoid-like structure on the spaces of Legendre fields. Note part (i) follows from the symmetry of $\nabla\circ$ alone - it does not require $\nabla e = 0$ - and so holds for generalized Frobenius manifolds.

\begin{proof}

\text{}

\begin{itemize}
\item[(i)] This follows immediately from the symmetry condition $(\nabla_X\circ)(e,e) = (\nabla_e\circ)(X,e)$ which gives $\nabla_Xe = X\circ \nabla_ee\,.$
\item[(ii)] Since $\tilde\nabla_XY = \delta^{-1}\circ \nabla_X(\delta\circ Y)$ it follows that $\tilde\nabla_X\delta^{-1} = \delta^{-1} \circ \nabla_X e$ and hence, from (i), that $\tilde\nabla_X\delta^{-1} = X \circ \delta^{-1} \circ \nabla_ee\,.$ Hence $\delta^{-1}$ is a Legendre field for $\tilde\nabla\,.$
\item[(iii)] The proof will be given in the context of Figure \ref{fig:transitive}, where $\widetilde\nabla$ is the connection for the generalized Frobenius manifold and, by definition of a Legendre field,
\begin{eqnarray*}
    \nabla^{(1)}_XY & = & \delta_1^{-1} \circ {\widetilde\nabla}_X (\delta_1\circ Y)\,,\\
    \nabla^{(2)}_XY & = & \delta_2^{-1} \circ {\widetilde\nabla}_X (\delta_2\circ Y)\,.
\end{eqnarray*}
On eliminating $\widetilde\nabla$ it follows that
\[
\nabla_X^{(2)}Y = \delta_2^{-1}\circ\delta_1\circ \nabla_X^{{(1)}} \left(\delta_1^{-1}\circ \delta_2\circ Y\right)\,.
\]
If the two Legendre fields are flat for $\tilde\nabla$, then $\nabla^{(i)}_Xe=0$. Hence $\delta_1^{-1}\circ\delta_2$ is a flat vector field for $\nabla^{(1)}$ and hence a Legendre field for this connection, and all fields arise in this way.
\end{itemize}
    
\end{proof}

\noindent The following is immediate (one just has to reverse some of Legendre fields in the above Proposition):
\begin{cor}\label{cor}
Let $\rm FM_1$ and $\rm FM_2$ be two Frobenius manifolds connected by a (flat, homogeneous) Legendre field $\delta\,:{\rm FM}_1 \xlongrightarrow{\delta} {\rm FM}_2\,.$ Let $\delta_1$ be a Legendre field for the connection $\rm FM_1\,.$ Then $\delta^{-1}\circ \delta_1$ is the Legendre field for $\rm FM_2\,.$
\end{cor}

\noindent The definition and properties of homogeneous Legendre field will be studied in Section \ref{section:homo} below, and it will immediately follow that if the two Legendre fields are homogeneous then the composite field $\delta_1^{-1}\circ\delta_2$ will also be homogeneous.

With a metric, to any vector field $V$ one may associate a form $\omega_V = \eta(V, -)$ (also denoted $V^\flat$). The following shows how the idea of a Legendre field, when coupled with a metric with the Frobenius property, is very natural.

\begin{lem}
    $(d\omega_V)(X,Y) = \eta( Y \circ \nabla_X V - X \circ \nabla_Y V,e)$ 
\end{lem}

\begin{proof}
    \begin{eqnarray*}
        (d\omega_V)(X,Y) & = & X[\omega_V(Y)] - Y[\omega_V(X)] - \omega_V([X,Y])\,,\\
        & = & X[\eta(V,Y)] - Y[\eta(V,X)] - \eta(V,[X,Y])\,,\\
        & = & \eta(\nabla_XV,Y) - \eta(\nabla_YV,X)\,,\\
        & = & \eta( Y \circ \nabla_X V - X \circ \nabla_Y V,e).
    \end{eqnarray*}
    This uses the fact that we have a torsion-free metric connection, and in the final line the Frobenius property of the metric has been used.
\end{proof}

\noindent Thus if $\delta$ is a Legendre field, then the associated form $\omega_\delta$ (which will be written as $\delta^\flat$) is closed and so $\delta^\flat = dh_\delta$ for some function $h_\delta$ which is defined up to an additive constant. We will call $h_\delta$ the potential for the Legendre field $\delta\,.$

\subsection{Homogeneous fields}\label{section:homo}
Let $\delta$ be a homogeneous Legendre field, that is
\[
\mathcal{L}_E \delta = \kappa\, \delta\,,
\]
where $\kappa$ is a constant (also called the degree of $\delta$). Since ${\mathcal{L}}_E\circ = \circ\,,$ then if Legendre fields $\delta_1\,,\delta_2$ have degrees $\kappa_1$ and $\kappa_2$ then $\delta_1\circ\delta_2$ has degree $\kappa_1+\kappa_2+1\,.$ We begin with some basic properties of such a field.

\begin{lem}\label{niceresult}
Let $h_\delta$ be a potential for the Legendre field $\delta$ (so $\delta^\flat = dh_\delta)$ and $\delta^{-1}$ its inverse\,. Then:
\begin{itemize}
    \item[(i)] $\mathcal{L}_E \delta^{-1} = (-2-\kappa) \, \delta^{-1}\,;$

 \vskip 2mm
    
    \item[(ii)] $(\delta^{-1})^\flat = \delta^\flat\,,$

    \item[(iii)] By suitable choice of the additive constant, the potential $h_\delta$ has the following homogeneity properties:
    \[\mathcal{L}_E h_\delta = 
    \begin{cases}
        (2-d+\kappa) \, h_\delta\,, & {\rm if~} 2-d+\kappa \neq 0\,,\\
       {\rm constant} \,, &{\rm if~} 2-d+\kappa =0\,,\\
    \end{cases}\]

    \vskip 2mm
        
    \item[(iv)] If $2-d+\kappa\neq 0$ then up to an inessential constant, $(2-d+\kappa) h_\delta = \eta(E,\delta)\,.$
\end{itemize}
\end{lem}

\begin{proof}

    \begin{itemize}
        \item[(i)] This follows from $\mathcal{L}_E \circ = \circ$ applied to $({\mathcal{L}}_E\circ)(\delta,\delta^{-1})=e\,.$ This also uses $[E,e]=-e$ which again follows from $({\mathcal{L}}_E\circ)(e,e)=e\,.$

        \item[(ii)] Care has to be taken in the interpretation of this. Since $\delta^{-1}$ is a Legendre field which acts in the opposite direction to $\delta$, the corresponding form is defined in terms of $\widetilde\eta$, so
                        \begin{eqnarray*}
                            (\delta^{-1})^\flat(X) & := & {\widetilde\eta}(\delta^{-1},X)\,,\\
                            & = & \eta(\delta,X)\,,\\
                            & = & (\delta)^\flat(X)\,.
                        \end{eqnarray*}
                        Note that this implies that $(\delta^{-1})^\flat$ and $\delta^\flat$ are both equal to the closed form $dh_\delta\,.$
        \item[(iii)] By definition,
                \begin{eqnarray*}
                    (\mathcal{L}_E \, dh_\delta)(X) & = & E[\eta(\delta,X)] - \eta(\delta,[E,X]) \,,\\
                    & = & ({\mathcal{L}}_E\eta)(\delta,X) + \eta(X,[E,\delta])\,,\\
                    & = & (2-d+\kappa) \, dh_\delta(X)\,.
                \end{eqnarray*}
            To find the homogeneity of $h_\delta$ one uses $[\mathcal{L},d]=0$ which gives $d ( \mathcal{L}_Eh_\delta - (2-d+\kappa) h_\delta)=0$ and hence $\mathcal{L}_Eh_\delta - (2-d+\kappa) h_\delta={\rm constant}\,.$ Since $h$ is only defined up to a constant, the result follows. 
                        
    \item[(iv)] Equation (\ref{skew}) with $Y=\delta$, namely $\eta(\mu(X),\delta) + \eta(X,\mu(\delta))=0\,,$ gives on using the torsion free condition and homogeneity,
    \[
    (2-d+\kappa) \eta(X,\delta) = \eta(\nabla_XE,\delta) + \eta (\nabla_E\delta,X)
    \]
    and using the metric connection condition on the first term on the right gives
    \begin{eqnarray*}
    (2-d+\kappa) \eta(X,\delta) & = & X(\eta(E,\delta)) + \eta(\nabla_E\delta,X)-\eta(E,\nabla X\delta)\,, \\
    & = & X(\eta(E,\delta)) + \eta(e, \underbrace{X\circ\nabla_E\delta - E\circ \nabla_X\delta}_{\text{$=0$ since $\delta$ Legendre}})\,.
    \end{eqnarray*}
    Hence
    \[
    (2-d+\kappa) dh_\delta(X) = \left(d\eta(E,\delta)\right)(X)
    \]
    so, up to constant which can be absorbed into $h_\delta$ (provided $2-d+\kappa\neq 0$), the result follows. 
    
    \end{itemize}
\end{proof}

\noindent A special case of  $(iv)$ with $\delta=e$ appeared in \cite{GF1}. In the following Proposition, Part $(i)$ was originally proved in \cite{SS} and the rest appeared in \cite{GF1}. Here, for completeness, we give a coordinate-free proof.

\begin{prop}\label{prop:spectrumhomogeneousLF}

Let $\delta$ be a homogeneous Legendre field with $\mathcal{L}_E \delta = \kappa \delta\,.$ Then:

\begin{itemize}
    \item[(i)] $\mathcal{L}_E \widetilde\eta = (2 - \tilde{d}) \, \widetilde{\eta}\,,\quad {\rm where~} \tilde{d}=d-2-2\kappa\,;$
    \vskip 2mm
    \item[(ii)] $\widetilde{\mu}(X) = \delta^{-1} \circ \mu ( \delta\circ X)\,;$
    \vskip 2mm
    \item[(iii)] $(\widetilde{\nabla}_X {\widetilde{\mu}})(Y) = \delta^{-1} \circ (\nabla_X\mu)(\delta\circ Y)\,;$
    \vskip 2mm
    \item[(iv)] The spectrum is invariant under the Legendre transformation defined by $\delta\,.$
\end{itemize}

\end{prop}

\begin{proof}

\phantom{text}

\begin{itemize}

\item[(i)] By direct calculation,
\begin{eqnarray*}
    (\mathcal{L}_E \tilde\eta)(X,Y)  =
    (\mathcal{L}_E \eta)(\delta\circ X,\delta \circ Y)&+&
    \eta([E,\delta\circ X] - \delta\circ [E,X],\delta\circ Y)\\
   &+&  \eta(\delta\circ X,[E,\delta\circ Y] - \delta\circ [E,Y])
\end{eqnarray*}
and from $(\mathcal{L}_E\circ)(\delta,X) = \delta\circ X$ it follows that $[E,\delta\circ X] - \delta\circ [E,X]=(1+\kappa) \delta\circ X$ and hence
\[
(\mathcal{L}_E \tilde\eta)(X,Y)=(2-d+2+2\kappa) \tilde{\eta}(X,Y)
\]
    \item[(ii)] By definition (\ref{endo}) of the endomorphism $\tilde\mu\,,$ 
    \begin{eqnarray*}
        \tilde{\mu}(X) & = & \frac{2-\tilde{d}}{2} X - \widetilde{\nabla}_X E\,,\\
        & = & \frac{2-\tilde{d}}{2} X - \widetilde{\nabla}_E X + [E,X]\,,\\
        & = & \delta^{-1} \circ \left\{ \frac{2-\tilde{d}}{2} \delta\circ X - \nabla_E(\delta\circ X) + \delta\circ[E,X] \right\}\,,\\
        & = & \delta^{-1} \circ \left\{ \frac{2-\tilde{d}}{2} \delta\circ X - \nabla_{\delta\circ X} E + \delta\circ[E,X] - [E,\delta \circ X]\right\}\,,\\
        & = &  \delta^{-1} \circ \mu(\delta\circ X)\,
        \end{eqnarray*}
        on using part (i) and $(\mathcal{L}_E\circ)(\delta,X) = \delta\circ X\,.$
\item[(iii)] Using the above result,
\begin{eqnarray*}
    (\tilde{\nabla}_X\tilde{\mu})(Y) & = & \tilde{\nabla}_X (\tilde\mu(Y)) - \tilde{\mu}(\tilde{\nabla}_XY)\,,\\
    & = & \delta^{-1}\circ \nabla_X(\mu(\delta\circ Y)) - \delta^{-1}\circ \mu(\nabla_X(\delta\circ Y))\,,\\
    & = & \delta^{-1} \circ (\nabla_X\mu)(\delta\circ Y)\,.
\end{eqnarray*}

\item[(iv)] The following results are standard:
\begin{eqnarray*}
\mathcal{L}_E\eta = (2-d)\eta & \Leftrightarrow & \mu {\rm{~skew~adjoint}}\,,\\
(\nabla_X\mu)(Y) & = & -  Hess^\nabla_{X,Y}(E)
\end{eqnarray*}
Thus if the metric is flat, and $E$ is linear in flat coordinates for this metric, then the Hessian is zero and hence $\nabla\mu=0\,.$ The above results show this argument can be reversed, and hence, since a Legendre transformation preserves any zero-curvature condition, the Euler field remains linear. A more detailed analysis shows that if
\[
E = \sum (d^i t^i + r_i)\frac{\partial~}{\partial t^i}
\]
where $d^i = 1 - \frac{d}{2} - \mu_i\,,$ then under a Legendre transformation the spectrum is invariant, i.e. $\tilde{\mu}_i = \mu_i\,.$

\end{itemize}
    
\end{proof}

\noindent How the two sets of flat coordinates for the two flat metrics $\eta$ and $\tilde{\eta}$ are related will be derived in Section \ref{coordinates}. With the above results the following is immediate:

\begin{cor}\label{Cor31}
    A homogeneous Legendre field maps a \gFM (with metric $\eta$ and connection $\nabla$) to a generalized Frobenius manifold (with metric $\widetilde\eta$ and connection $\widetilde\nabla$) where:
    \begin{eqnarray*}
        \widetilde{\eta}(X,Y) & = & \eta(\delta\circ X,\delta\circ Y)\,;\\
        \widetilde{\nabla}_XY & = & \delta^{-1} \circ \nabla_X(\delta\circ Y)\,;\\
        \tilde{d} & = & d - 2 - 2 \kappa\,;\\
        \tilde{\mu}_i & = & \mu_i\,
    \end{eqnarray*}
with $E$ and $e$ remaining unchanged. Moreover, if $\delta$ is flat (with respect to $\nabla$) and homogeneous then the field maps a \gFM to a Frobenius manifold whose charge $\tilde{d}$ depends on the spectrum of the generalized Frobenius manifold.
\end{cor}

\begin{proof}
The only part to prove is the last. If $\delta$ is flat with respect to $\nabla$ and homogeneous, then
\[
    \delta = \frac{\partial~~}{\partial t^\bullet}
\]
for some coordinate denoted by $t^\bullet$ (if the elements of spectrum are not distinct one can redefine the flat coordinates appropriately) and this has homogeneity $\kappa =-d^\bullet$ and in terms of the spectrum (\ref{spectrum})
\[
\kappa = -1 + \mu_\bullet + \frac{d}{2}\,.
\]
Since $\nabla_X \delta=0$ then $\tilde{\nabla}_X e = \delta^{-1}\circ \nabla_X \delta= 0\,$
and hence the new structure is a Frobenius manifold, not a generalized Frobenius manifold. The last two conditions then imply that 
\[
\tilde{d} = -2 \mu_\bullet\,,
\]
so the charge of the Frobenius manifold obtained by this construction is independent of the charge of the generalized Frobenius manifold - only on the spectrum. It also follows that the coefficients in the Euler field transform as
\[
\tilde{d}^i = d^i - (1+\kappa)\,.
\]
\end{proof}

\begin{ex}\label{ex3.2}
    Consider the generalized Frobenius manifold with Euler field and charge
    \begin{eqnarray*}
        E & = & \left(1-\frac{1}{N}\right) \tilde{t}^1\ \frac{\partial~}{\partial\tilde{t}^1 } + \frac{1}{N} \tilde{t}^2\ \frac{\partial~}{\partial\tilde{t}^2 }\,,\\
        \tilde{d} & = & 1\,
    \end{eqnarray*}
    (the two dimensional examples in \cites{orbit1,orbit2} fall into this class). There are two choices for flat, homogeneous Legendre fields - the two coordinate vector fields - and the above formulae gives (on relabeling the coordinates) the following:

    \begin{itemize}
        \item[(a)] if
        $\displaystyle{\delta_1 = \frac{\partial~}{\partial\tilde{t}^1 }}\,$  then
        \begin{eqnarray*}
        E & = &  {t}^1 \frac{\partial~}{\partial{t}^1 } + \frac{2}{N} {t}^2 \frac{\partial~}{\partial{t}^2 }\,,\\
        {d} & = & 1-\frac{2}{N}\,;
        \end{eqnarray*}
        \item[(b)] if
        $\displaystyle{\delta_2 = \frac{\partial~}{\partial\tilde{t}^2 }}\,$  then
        \begin{eqnarray*}
        E & = &  {t}^1 \frac{\partial~}{\partial{t}^1 } + \left(2-\frac{2}{N}\right) {t}^2 \frac{\partial~}{\partial{t}^2 }\,,\\
        {d} & = & \frac{2}{N}-1\,;
        \end{eqnarray*}
    \end{itemize}
This illustrates (see Figure \ref{fig:transitive}) the fact that a generalized Frobenius manifold is \lq above\rq~different Frobenius manifolds that are themselves related by a flat, homogeneous, Legendre transformation.

\end{ex}

\subsection{Twisted Legendre transformation and almost-duality}

Given a Frobenius manifold with metric $\eta$ and multiplication $\circ$ one may construct an associated manifold - an almost-dual Frobenius manifold - with new metric and multiplication
\begin{eqnarray*}
g(X,Y) & = & \eta(E^{-1}\circ X,Y)\,,\\
X\star Y & = & E^{-1}\circ X\circ Y
\end{eqnarray*}
where $E$ is the Euler vector field and $E^{-1}\circ E=e$ (we assume we are working on an open set in the manifold in which $E$ is invertible). It was shown in \cites{Dubrovin1996,Dub04} that $g$ is flat and compatible with the new product and, moreover, that a dual solution of the WDVV equations (written in the flat coordinates of the metric $g$) may be constructed. While the metrics $\eta$ and $g$ are related by the simple formula above, the corresponding connections, $\nabla$ and $\nabla^\star$, are related by a more complicated formula (see, for example, \cites{HertlingBook,davidstrDubduality})
\begin{equation}\label{dualconnection}
\nabla^\star_XY = E\circ\nabla_X(E^{-1}\circ Y) - X\circ\nabla_{E^{-1}\circ Y} E + \frac{1}{2} (3-d) E^{-1}\circ X\circ Y\,.
\end{equation}
These proofs do not use the condition $\nabla_Xe=0$ and hence still hold for a \gFM - one can construct an almost-dual \gFM via the same construction.

Since we have an almost-dual solution of the WDVV equations, we automatically have a flat pencil of connections
\[
{}^\lambda\nabla^\star_XY = \nabla^\star_XY - \lambda X\star Y
\]
and hence the above theory of generalized Legendre transformations may be applied directly to solutions by replacing $\eta$ with $g$ and $\nabla$ with $\nabla^\star\,.$ 

Thus one has, for a homogeneous Legendre field $\delta\,,$ the following schematic diagram:

\[
\begin{array}{cccc}
&{\rm gFM}& \overset{\delta}{\xrightarrow{\hspace{1cm}}}& \widetilde{{\rm gFM}}\\
& & \\[-5pt]
\text{\footnotesize{almost duality}}&\big\downarrow & &\big\downarrow  \\
& & \\[-5pt]
&{\rm gFM}^\star &&~~~~\widetilde{\rm gFM}^\star\,.
\end{array}
\]

\bigskip

\noindent The next results construct a Legendre field $\delta^\star$ that makes the above diagram commute:
\[
{\rm gFM}^\star ~\overset{\delta^\star}{\xrightarrow{\hspace{1cm}}} ~\widetilde{\rm gFM}^\star\,.
\]
This field - called a twisted Legendre field - was first studied in \cite{RS}. However the full geometric properties were not fully described there.

\begin{thm}\label{thm:almostdual}
The following diagram commutes:
\[
\begin{array}{ccc}
(\eta\,,\nabla\,,\circ)& \overset{\delta}{\xrightarrow{\hspace{1cm}}}&{({\tilde\eta}\,,{\tilde\nabla}\,,\circ)}\\
& & \\
\big\downarrow & &\big\downarrow  \\[-10pt]
& & \\
(g\,,\nabla^\star\,,\star) & \overset{\delta^\star} {\xrightarrow{\hspace{1cm}}}& ~~({\tilde g}\,,{\tilde \nabla}^\star,\star)\,,
\end{array}
\]
where $\delta^\star = E\circ\delta$ and $\delta$ is a homogeneous Legendre field with $\mathcal{L}_E\delta = \kappa \delta\,.$ In particular:
\begin{itemize}
\item[(i)] ${\tilde g}(X,Y)= g ( {\delta^\star}\star X, {\delta^\star}\star Y)$\,;
\vskip 2mm
\item[(ii)] $\mathcal{L}_E{\delta^\star} = (1+\kappa) \delta^\star\,.$
\vskip 2mm
\item[(iii)] $\delta^\star$ is a Legendre field for the connection $\nabla^\star\,;$
\vskip 2mm
\item[(iv)] ${\tilde\nabla}^\star_XY = {\delta^\star}^{-1} \star \nabla_X^\star ( {\delta^\star}\star Y)$\,;
\vskip 2mm
\item[(v)] $(\delta^\star)^\flat = \delta^\flat\,.$
\end{itemize}

\end{thm}

\begin{proof} 

\phantom{text}

\begin{itemize}
\item[(i)] By definition
\[
{\tilde g}(X,Y) = {\tilde \eta}(E^{-1}\circ X,Y) 
= \eta(\delta\circ E^{-1}\circ X,\delta\circ Y)
=g(\delta\circ X,\delta\circ Y)
= g({\delta^\star}\star X, {\delta^\star}\star Y)\,.
\]
since $\delta\circ X = {\delta^\star} \star X\,.$ 

\item[(ii)] This follows immediately from $(\mathcal{L}_E\circ)(E,\delta) = E\circ\delta = \delta^\star\,.$

\item[(iii)] From (\ref{dualconnection}),
\[
\nabla_X^\star(\delta\circ E) = E\circ \nabla_X\delta - X \circ \nabla_\delta E + \frac{1}{2}(3-d) X\circ \delta\,.
\]
Hence $Y\star\nabla_X^\star(\delta\circ E)$ is symmetric in $X$ and $Y$ since $\delta$ is a Legendre field. Furthermore
\begin{eqnarray*}
    \nabla_X^\star \delta^\star & = & E\circ \nabla_X\delta - X \circ \nabla_E \delta + X\circ [E,\delta] + \frac{1}{2}(3-d) X\circ \delta\,\\
    & = & \left( \frac{1}{2}(3-d) + \kappa\right) X\circ\delta\,\\
    & = & \frac{1}{2} \left(1-\tilde{d}\right) X \circ \delta
\end{eqnarray*}
since $\delta$ is a Legendre field. It then follows that $Y\star\nabla_X^\star\delta^\star$ is symmetric in $X$ and $Y\,,$ and hence $\delta^\star$ is a Legendre field for the connection $\nabla^\star\,.$

\item[(iv)] Starting with the connection at the lower right hand corner and unpacking definitions yields
\begin{eqnarray*}
{\tilde\nabla}^\star_XY &=& E\circ{\tilde\nabla}_X(E^{-1}\circ Y) - ({\tilde\nabla}_{E^{-1}}\circ Y)\circ X + \frac{1}{2} (3-{\tilde d}) E^{-1}\circ X\circ Y\,,\\
&=&E\circ\left\{
\delta^{-1} \circ \nabla_X(\delta \circ E^{-1} \circ Y) \right\}-
\left\{\delta^{-1} \circ \nabla_{E^{-1}\circ Y} (\delta\circ E)\right\}\circ X \\
&&+ \frac{1}{2} (3-{\tilde d}) E^{-1}\circ X\circ Y\,,\\
&=&\delta^{-1}\circ \nabla^\star_X(\delta\circ Y) + \delta^{-1} \circ\left\{ \nabla_{E^{-1}\circ \delta\circ Y}E - \nabla_{E^{-1}\circ Y}(\delta\circ E)\right\}\circ X\\
&&+ \frac{1}{2} (d-{\tilde d}) E^{-1}\circ X\circ Y\,.
\end{eqnarray*}
On using the following formula derived in \cite{DavidStrDdualityconnection}

\begin{eqnarray*}
\nabla_{E^{-1}\circ \delta\circ Y}E - \nabla_{E^{-1}\circ Y}(\delta\circ E) &=& \left\{[\delta,E]+[E,e]\circ\delta\right\}\circ Y\circ E^{-1}\,,\\
&=& -(1+\kappa)\,\delta\circ Y\circ E^{-1}\,.
\end{eqnarray*}
Hence
\[
{\tilde\nabla}^\star_XY = \delta^{-1} \circ \nabla^\star_X(\delta \circ Y) + \left\{ -\kappa-1+\frac{1}{2}(d-{\tilde d})\right\}\,E^{-1}\circ X\circ Y\,
\]
and the coefficient of the second terms vanishes. Hence $(\widetilde\nabla)^\star= \widetilde{(\nabla^\star)}\,.$

\item[(v)] Again, care has to be taken with the definition of $(\delta^\star)^\flat\,,$ which uses the metric $g\,.$ Thus
\[
(\delta^\star)^\flat(X) = g(\delta^\star,X) = \eta(E^{-1}\circ E\circ \delta,X) = \delta^\flat(X)\,.
\]
Hence $(\delta^\star)^\flat = dh = \delta^\flat\,,$ with the interpretation that different metrics are used in defining the 1-forms.

\end{itemize}

\end{proof}

\noindent Note that the statement of Theorem \ref{thm:almostdual}(v) extends: on defining $(\delta^\star)^{-1}$ using the dual multiplication $\star$ gives $(\delta^\star)^{-1} = E\circ\delta^{-1}$ and it follows from this that $\left((\delta^\star)^{-1}\right)^\flat$ (defined using the metric $\tilde{g}$) is equal to $\delta^\flat$ and hence is also equal to $dh\,.$ Thus all the connecting Legendre fields, in either direction, have the same potential $h\,.$

\subsection{Coordinate representations}\label{coordinates}

So far we have used a coordinate-free approach. However, for examples one requires explicit computational formulae in the flat coordinates of the various metrics. This is provided by the following Theorem that extends the results in \cite{SS}.

\begin{thm}\label{symmetry}
Let $\{t^\alpha\}$ and $\{ {\tilde t}^\alpha\}$ be the flat coordinates of the metrics $\eta$ and $\tilde\eta$ respectively. Then, given a \gFM with prepotential $F$ and a homogeneous Legendre field $\delta$ with potential $h$, the new flat coordinates, metric coefficients and prepotential for the Legendre-transformed \gFM are given by:
\begin{itemize}
\item[(i)] $\displaystyle{{\tilde\eta}_{ij} = \eta_{ij}}\,;$
\vskip 2mm
\item[(ii)] If $\tilde{d}^i\neq 0$ then $\tilde{d}^i \,{\tilde{t}^i} = g^{ij}(t)\,\partial_j h(t)\,;$
\vskip 2mm
\item[(iii)] 

$\displaystyle{\frac{\partial^2 \widetilde{F}}{\partial\tilde{t}^i \partial\tilde{t}^j} =  \frac{\partial^2 {F}}{\partial{t}^i \partial{t}^j} }\,.$

\end{itemize}
\end{thm}

\begin{proof}

\begin{itemize}
\item[(i)] Since $\tilde\eta(X,Y) =\eta(\delta\circ X,\delta \circ Y)$ it follows that, up to an inessential linear transformation,
\begin{equation}
\frac{\partial~}{\partial t^i} = \delta \circ \frac{\partial~}{\partial {\tilde t}^i}
\label{CoV}
\end{equation}
\noindent and hence that
\begin{eqnarray*}
{\tilde\eta}_{ij} & = & {\tilde\eta}\left( \frac{\partial~}{\partial {\tilde t}^i},\frac{\partial~}{\partial {\tilde t}^j}\right)\,,\\
&=&  {\eta}\left( \frac{\partial~}{\partial {t}^i},\frac{\partial~}{\partial {t}^j}\right)\,,\\
&=& \eta_{ij}\,.
\end{eqnarray*}

\item[(ii)] From (\ref{CoV}) it follows that
\begin{eqnarray*}
\frac{\partial~}{\partial t^i} & = & \frac{\partial t^j}{\partial {\tilde t}^i} \left( \delta \circ 
\frac{\partial~}{\partial {t}^j}\right)\,,\\
&=&\frac{\partial t^j}{\partial {\tilde t}^i} \, \delta^k \,c_{k}^r\, \frac{\partial~}{\partial t^r}\,.\\
\end{eqnarray*}
Hence
\begin{equation}
\frac{\partial {\tilde t}^i }{\partial t^j} = \delta^k c^i_{kj}
\label{changeofvariable}
\end{equation}
Since $\delta^k=\eta^{ks}\partial_s h(t)$ contracting with the components of the Euler field $E^j$
gives, on using the definition of the intersection form, 
\begin{eqnarray*}
E({\tilde{t}}^i) & = &  c^{ij}_k(t)  \, E^k(t) \,\partial_j h(t)\,,\\
& = & g^{ij}(t)\,\partial_j h(t)
\end{eqnarray*}
and since $E({\tilde t}^i) = \tilde{d}^i \, {\tilde{t}^i}$ the result follows. Note, this construction only gives a closed form expression for $\tilde{t}^i$ if $\tilde{d}^i\neq 0
\,.$

\item[(iii)] To prove the final formula,
\begin{eqnarray*}
\frac{\partial {\tilde F}_{ij}}{\partial {\tilde t}^k} = {\tilde c} \left(
\frac{\partial~}{\partial {\tilde t}^i},
\frac{\partial~}{\partial {\tilde t}^j},
\frac{\partial~}{\partial {\tilde t}^k}\right) & = & {\tilde{\eta}}
\left(\frac{\partial~}{\partial {\tilde t}^i}\circ\frac{\partial~}{\partial {\tilde t}^j},\frac{\partial~}{\partial {\tilde t}^k}\right)\,,\\
&=&
\eta\left(
\frac{\partial~}{\partial {t}^i}\circ\frac{\partial~}{\partial {t}^j},\frac{\partial~}{\partial {t}^r}\right) \frac{\partial t^r}{\partial {\tilde t}^k}\,,\\
&=&
\frac{\partial t^r}{\partial {\tilde t}^k}
c\left(
\frac{\partial~}{\partial {t}^i},
\frac{\partial~}{\partial {t}^j},
\frac{\partial~}{\partial {t}^r}\right)
= \frac{\partial F_{ij}}{\partial {\tilde t}^k} \,.
\end{eqnarray*}
Hence, on integrating (and recalling that the prepotentials are only defined up to quadratic terms in the appropriate flat coordinates) part (c) follows.
\end{itemize}
\end{proof}

\noindent The appearance of the intersection form is perhaps surprising. Note that, for Frobenius manifolds, $g^{ni}=E^i$ and hence, up to multiplicative constants, $(ii)$ gives $\tilde{t}^n=h(t)\,.$ This result also provides computational formula to perform Legendre transformations. The construction, however, require an inversion to get $t^i({\bf \tilde{t}})$ which, even in the simplest low-dimensional examples, is hard.

\subsection{Canonical Legendre transformations}

The fact that $\delta^{\star\flat} = \delta^\flat = dh$ has not so far been used. This is now used to provide a canonical class of Legendre fields.

\begin{thm}
Consider a \FM with metric $\eta$ and flat coordinates 
$\{t^i\}$ and intersection form $g$ with flat coordinates $\{p^i\}\,.$ A class of Legendre fields $\delta$ may be constructed from the potential
\[
h({\bf t}) = \alpha_i \,p^i({\bf t})
\]
where $\alpha_1\,,\dots\,,\alpha_n$ are constants, and this gives the Legendre field $\delta = dh^\sharp\,.$ Moreover, the resulting \gFM has $\tilde{d}=1\,.$ Conversely, all generalized Frobenius manifolds with charge $\tilde{d}=1$ arise from this construction. 
\end{thm}

\begin{proof}
    A simple class of Legendre fields for the almost-dual connection $\nabla^\star$ are just linear combinations of the flat vector fields for the metric (or intersection form) $g\,,$ namely
    \[
    \delta^\star = \alpha_i \frac{\partial~~}{\partial p^i} 
    \]
    where the $\alpha_i$ are constants, and this has the potential $h=\alpha_i p^i\,.$ Since, from Theorem \ref{thm:almostdual}(v) one has $\delta^{\star\flat} = \delta^\flat = dh$, the first result follows. The second part follows from the details of the proof of Theorem \ref{thm:almostdual}(iii) where it was shown that
    \[
    \nabla_X^\star \delta^\star  = \frac{1}{2} \left(1-\tilde{d}\right) X \circ \delta\,
    \]
and hence $\delta^\star$ is flat for $\nabla^\star$ if and only if $\tilde{d}=1\,.$ Using this, the Legendre field $\delta^\star$ must have $\tilde{d}=1$, and conversely if the Legendre field generates a generalized Frobenius manifolds with $\tilde{d}=1$ then the twisted field $\delta^\star$ must be flat for the dual connection.
    
%
%
%

\end{proof}

\medskip

\noindent Since we require $\delta$ to be invertible the constants $\alpha^i$ cannot all be zero, and as the scale does not play a role, the free data is a point in $\mathbb{P}^{n-1}\,.$ It is important to note that the flat coordinates for the intersection form are not uniquely defined and so this construction does not single out distinguished points in this space. This construction, though, circumvents the need to solve the Legendre field equations. One does, though, have to construct the flat coordinates of $g(t)$ in terms of the flat coordinates of $\eta\,.$ This involves solving the Gauss-Manin equations $\nabla^\star \,dp(t)=0$ to obtain $p^i=p^i({\bf t})\,.$ Solutions of these Gauss-Manin equations are multi-valued, branching over the discriminant
\[
\Sigma = \{ {\bf{t}}\,:\, \det\left(g^{ij}(t)\right)=0\}
\]
and this multi-valuedness is best described in terms of the monodromy group 
$\mu\left[ \pi_1(M\backslash \Sigma) \right]$ of the Frobenius manifold \cite{Dubrovin1996}. We will also see these monodromy groups when we directly solve the Legendre field equations in Section \ref{classification}. Also, in constructions such as the Saito construction for the Coxeter group orbit space $\mathbb{C}^n/W\,,$ one obtains the flat coordinates on the orbit space as Saito polynomials $t^i\in \mathbb{C}[{\bf p}]^W$ and hence to construct $h$ one has to invert these polynomials. 

This result singles out \gFMs with $\tilde{d}=1$ as, in some sense, distinguished. This motivates the following:

\begin{defn}
A canonical \gFM is a generalized Frobenius manifold with charge $\tilde{d}=1\,.$
\end{defn}

From the transitive action of Legendre transformations any \gFM may be transformed into a canonical generalized Frobenius manifold, and the above Theorem gives a construction that takes one from an arbitrary \FM to a canonical generalized Frobenius manifold. Another reason why having $\tilde{d}=1$ is special will be given in the next section. 

\begin{ex}
    Consider the Frobenius manifold $FM_1$ (see Table 1) on (for $N\in \mathbb{N}_{\geq 3}$) the Coxeter group orbit space $\mathbb{C}^2/I_2(N)\,.$ The flat coordinates for the intersection form are \cite{Dubrovin1996}
    \begin{eqnarray*}
        t^1 & = & \frac{1}{N^\frac{N+2}{2}} \,(z^N+{\overline{z}^N})\,,\\
        t^2 & = & \frac{1}{N}\, z {\overline{z}}
    \end{eqnarray*}
    and $g=2\, dz \, d\overline{z}$ (the normalization are for later use). These can easily be inverted to give the canonical Legendre field potential
    \[
    h(t^1,t^2) = \alpha_1 \left(
    t^1 + \sqrt{(t^1)^2 - \frac{4}{N^2} (t^2)^N}\right)^{\frac{1}{N} } + \alpha_2 \left(
    t^1 - \sqrt{(t^1)^2 - \frac{4}{N^2} (t^2)^N}\right)^{\frac{1}{N} }\,.
    \]
    Note the branching along the discriminant $\frac{4}{N^2} (t^2)^N- (t^1)^2=0$ which is where $det(g^{ij})=0\,.$ This may be rewritten as $h(t^1,t^2) = (t^2)^{1/2} \, y(x)$ where
    \[
    y(x) = \alpha_1 \left(
    x + \sqrt{x^2 - 1}\right)^{\frac{1}{N} } +
    \alpha_2 \left(
    x - \sqrt{x^2-1}\right)^{\frac{1}{N} } 
    \]
    where $x$ is the degree zero variable $x=\frac{N t^1} {2 (t^2)^{N/2}}\,.$ This is an algebraic function and hence satisfies an algebraic relation $P^N_{[\alpha_1,\alpha_2]}(x,y) = 0$ and one can see that generically this is of degree $2N$ in $y\,.$ 

The local monodromy matrices at the points $x=+1$ and $x=-1$ are just
    \[
    M_1 = \left(\begin{array}{cc} 0 & 1 \\ 1 & 0 \end{array}\right)\,, \qquad M_{-1} = \left(\begin{array}{cc} 0 & 1 \\ 1 & 0 \end{array}\right)
    \]
    and the eigenvectors of these simultaneously diagonalizable matrices single out the two subspaces $\alpha_1-\alpha_2=0$ and $\alpha_1+\alpha_2=0$ as special (and as eigenspaces, they are basis independent). While the local monodromy group is $\mathbb{Z}_2$, its representation degenerates to the trivial representation in the first case, and to minus the identity in the second. We can also identify these special cases via the geometry of the curve $P^N_{[\alpha_1,\alpha_2]}(x,y) = 0$ (for a discussion of the full monodromy group, see Example \ref{FM1example}).

\vskip 2mm

\noindent{\underline{Case 1: $\alpha_1-\alpha_2=0$}}
On setting, without loss of generality, $\alpha_1=1/2\,,$ then
if $x=\cos N\theta$ then $y=\cos\theta$. The defining formula for the Chebyshev polynomials of the first kind, $T_N(\cos\theta) = \cos\,N\theta\,,$ then gives the simple algebraic curve
\[
T_N(y) = x\,.
\]
which is of lower degree than the generic case.

\vskip 2mm

 \noindent{\underline{Case 2: $\alpha_1=-\alpha_2$}}
 Similarly, if $x=\sin N\theta$ then $y=\cos\theta$ then using the formula for Chebyshev polynomials of the second kind,
\[
U_{N-1}(\cos\theta) = 
    \displaystyle{\frac{\sin N\theta}{\sin\theta}} \,,
\]
and hence the identity $T_N^2(y) +(1-y^2) U_{N-1}^2(y)=1\,,$
one obtains the algebraic curve
\[
T_N^2(y) + x^2 = 1\,.
\]
Since $T_N(-y)=(-1)^N \,T_N(y)$ this curve depends on $x^2$ and $y^2$, so in terms of a covering map (replacing $y^2$ by $y$ and $x^2$ by $x$) one obtains a curve
\[
\widetilde{T}_N(\tilde{y}) = \tilde{x}\,,
\]
where $\widetilde{T}_N$ is a polynomial of degree $N$ (this construction does introduce branch points). Alternative one has a curve with a higher degree of symmetry than the generic case.

\vskip 2mm

\noindent Thus in the space of Legendre fields we have distinguished points and hence distinguished generalized Frobenius manifolds. There is a third class of special points not detected by local monodromy.

\vskip 2mm

\noindent{\underline{Case 3: $\alpha_1=0$ or $\alpha_2=0$ }}
The above formula immediately gives
\[
y^{2N} - 2 x y^N + 1 =0\,.
\]
\end{ex}
\noindent Again, in terms of covering, this reduces to a quadratic curve.

\vskip 2mm

Thus for special values of the constants the degree of the algebraic curve drops from its generic value or the curve becomes particularly symmetric. We see how these are related to certain explicit generalized Frobenius manifolds constructed in \cites{orbit1,orbit2}. The Frobenius manifolds $FM_2$ and $FM_3$ work in a similar fashion to $FM_1$ and we omit the details.

\begin{ex}\label{deltaP1}
    Consider the Frobenius manifold $FM_4$ (see Table 1) defined on the extended affine Weyl orbit space $\mathbb{C}^2/A^{(1)}_1\,.$ 
The transformation between the flat coordinates of the Saito and intersection form is well-known \cite{Dubrovin1996}, and the above construction results in the Legendre field
\[
\delta = \alpha_1 \frac{\partial~~}{\partial t^1} + \alpha_2
\Delta^{-1/2}\left\{
t^1\frac{\partial~~}{\partial t^1} - 2 
\frac{\partial~~}{\partial t^2}
\right\}
\]
where $\Delta = (t^1)^2 - 4 e^{t^2}\,.$ This is particularly simple since $d=1$ as well as $\tilde{d}=1\,.$ 
\end{ex}

\begin{figure}[H]
\begin{tikzpicture}
\node at (-2,1.5) {$\left\{\rule{0pt}{1.5cm}\right.$};
\draw[->] (-0.2,0.1) -- (-0.2,5.5);
\draw (-0.2,3.95) -- (1.0,5.15) -- (4.5,5.15) -- (3.3,3.95) -- cycle;
\draw (-0.2,1.5) -- (1.0,2.7) -- (4.5,2.7) -- (3.3,1.5) -- cycle;
\node[] at (6.0,5.5) {\footnotesize Canonical generalised};
\node[] at (6.0,5.0) {\footnotesize Frobenius manifolds};
\node[] at (-1.0,1.5) {${\tilde{d}}=-2\mu_i$};
\draw (-0.2,2.5) -- (0.2,2.9);
\draw (-0.2,2.5) -- (0.2,2.5);
\node[] at (-1.0,2.5) {${\tilde{d}}=-2\mu_1$};
\draw (-0.2,0.5) -- (0.2,0.9);
\draw (-0.2,0.5) -- (0.2,0.5);
\node[] at (-1.0,0.5) {${\tilde{d}}=-2\mu_n$};
\draw[thick,->] (0,-0.3) parabola (1.5,2); 
\node[] at (-0.53,-0.3) {FM $\cross$};
\node[] at (1.88,2.2) {$\cross$ FM};
\node[] at (1.8,1) {$\delta=\partial_{t^i}$};
\node[] at (-0.7,4.0) {${\tilde{d}}=1$};
\node[] at (-3.5,1.6) {$n$-distinguished};
\node[] at (-3.3,1.25) {values};
\end{tikzpicture}
\caption{Flat homogeneous Legendre fields (coordinate vector fields) lie in distinguished subspaces in the space of all Legendre fields.}
\label{fig:subspaces}
\end{figure}
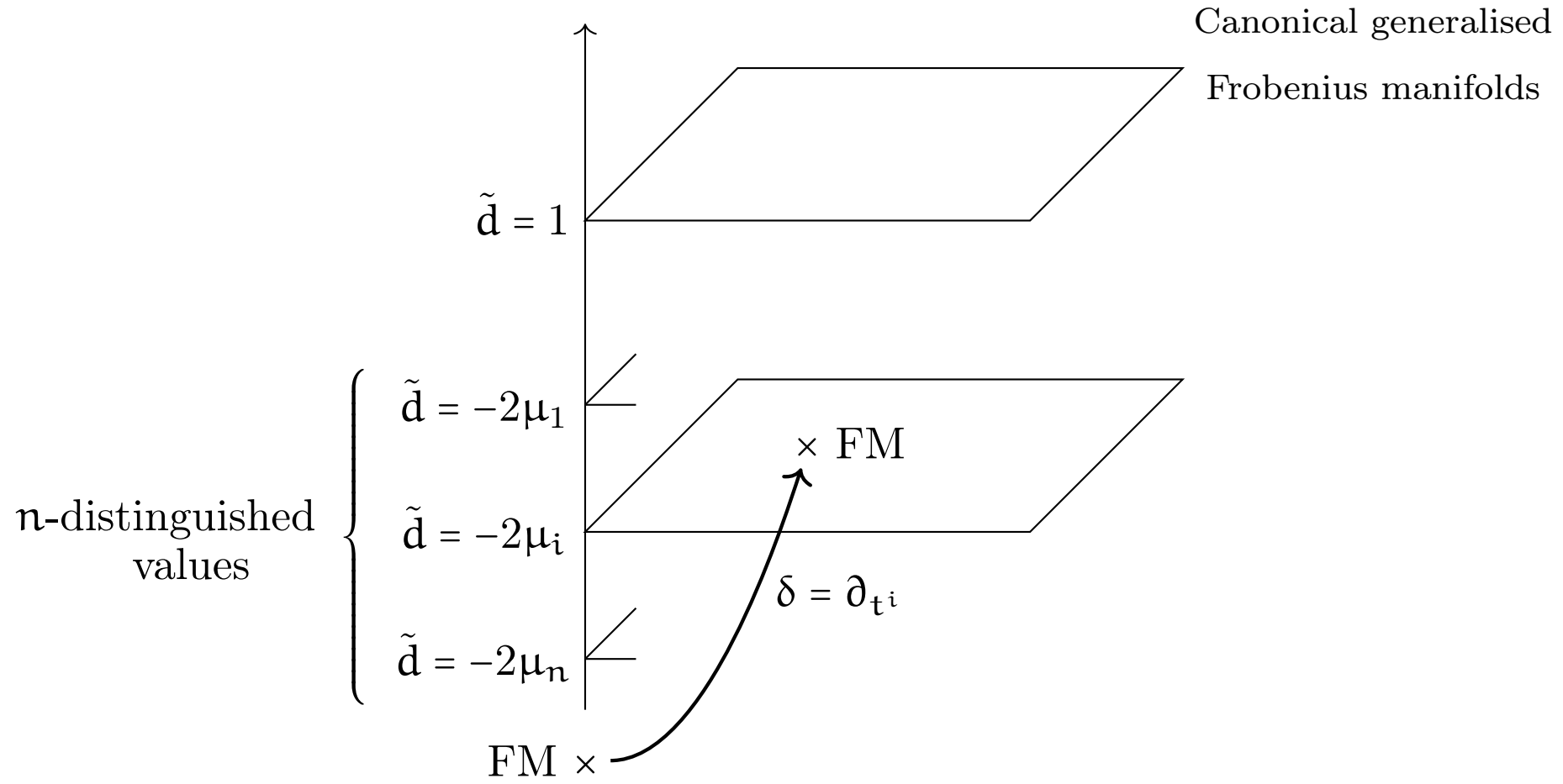

\noindent In this example the space of Legendre fields includes the flat homogeneous vector field $\partial_{t^1}\,.$ In general, the space of Legendre fields constructed from a specific Frobenius manifold will always include the $n$-coordinate vector fields $\partial_{t^i}\,.$ Since these have homogeneity $-d_i\,,$ it is easy to find the corresponding values of $\tilde{d}$ and these are given in terms of the spectrum of the underlying Frobenius manifold,
\begin{eqnarray*}
\tilde{d} & = & d - 2 + 2 d_i\,,\\
& =&  d-2+2 \left[1-\left(\mu_i + \frac{d}{2}\right)\right]\,,\\
& = & -2 \mu_i\,.
\end{eqnarray*}
This gives $n$-distinguished (and, if there is no degeneracies amongst the degrees, only $n$) subspaces in the space of Legendre fields (see Figure \ref{fig:subspaces}) - these containing the flat, homogeneous coordinate vector fields that map to Frobenius manifolds. If $-\frac{1}{2}$ belongs to the spectrum (as it does for the extended affine Weyl Frobenius manifold $FM_4$ used in Example \ref{deltaP1}), then vectors in this subspace will in general generate canonical generalised Frobenius manifolds. 

We end this section with the following, which brings together various results:

\begin{thm}\label{mapbetweenspaces}

\phantom{text}

\begin{itemize}
    \item[(a)]
    Consider two Frobenius manifolds ${\rm FM}_{1,2}$ with corresponding homogeneous Legendre fields $\delta_{1,2}$ which are related by a (flat, homogeneous) Legendre field
    \[
    {\rm  FM}_1 \xlongrightarrow{\hspace{3mm}{\frac{\partial~~}{\partial t^k}} \hspace{3mm}} 
    {\rm FM}_2\,.
    \]
    This then lifts to a homogeneous Legendre field $\delta$ 
    between the corresponding generalized Frobenius manifolds defined by
    \[
    \delta = \delta_1^{-1} \circ\left({\frac{\partial~~}{\partial t^k} } \right)\circ \delta_2\,.
    \]
    If these two generalized Frobenius manifolds have the same charge then
    \[
    \mathcal{L}_E \delta = - \delta\,.
    \]
    and this provides 
    isomorphism on the corresponding space of Legendre fields
    \[
    \delta:\mathbb{P}_{d,\tilde{d}}^{n-1} \xlongrightarrow{\hspace{3mm} \hspace{3mm}} \mathbb{P}_{d-2+2d_k,\tilde{d}}^{n-1}
    \]
 where $d$ is the charge of $FM_1$ and $\tilde{d}$ the charge of the generalized Frobenius manifold constructed from $\delta_1\,.$

 \item[(b)] Consider two generalized Frobenius manifolds ${\rm gFM_{1,2}}$ with the same spectrum and two (flat homogeneous) Legendre fields mapping to the Frobenius manifolds ${\rm FM}_{1,2}\,.$ Then if a Legendre transformation is to exists between the Frobenius manifolds it is uniquely determined by the spectrum.

 \end{itemize}
    
\end{thm}

\begin{oss}
Before proving this theorem it is important to point out the following: it is only in one and two dimensions that Frobenius manifolds are uniquely characterized by their charge. In higher dimensions one has other data to consider, such as the monodromy data, Stokes and central connection matrices. While a Legendre transformation preserves such data (see \cite{yang}) one cannot deduce that two Frobenius manifolds with the same charge and spectrum are related by a Legendre transformation. This is already apparent in three dimensions, where the WDVV equation reduce, for a specified charge and spectrum, to a third ordinary differential equation \cite{Dubrovin1996}. With different boundary/initial data one obtains different solutions and hence different manifolds. This explains the wording in the part (b) of the Theorem. 
\end{oss}

\begin{proof} 

\begin{itemize}

    \item[(a)]

Schematically we have the following:
\begin{figure}[H]
\begin{tikzpicture}
\draw[->] (-0.2,0.1) -- (-0.2,2.8);
\draw (-0.2,1.5) -- (1.0,2.5) -- (3.5,2.5) -- (2.3,1.5) -- cycle;
\draw[thick,->] (0,0) parabola (1.5,1.8); 
\node[] at (-0.6,0.02) {$\rm FM_1$ $\cross$};
\node[] at (1.55,2.0) {$\cross$};
\draw[->,dashdotted] (1.8,2) -- (8.5,2);
\node[] at (1.4,1) {$\delta_1$};
\node[] at (-0.4,1.55) {${\tilde{d}}$};

\draw[->] (7,0.1) -- (7,2.8);
\draw (7,1.5) -- (8.2,2.5) -- (10.7,2.5) -- (9.5,1.5) -- cycle;
\draw[thick,->] (7.2,0) parabola (8.7,1.8); 

\node[] at (6.6,0.02) {$\rm FM_2$ $\cross$};
\node[] at (8.75,2.0) {$\cross$};
\node[] at (8.6,1) {$\delta_2$};
\node[] at (6.8,1.55) {${\widehat{d}}$};

\draw[thick, ->] (1,0) -- (5.8,0);
\node[] at (3.6,0.5) {$ \frac{\partial~~}{\partial t^k} $};

\end{tikzpicture}
\label{fig:inducedmap}
\end{figure}
\noindent From Proposition \ref{transitive} it follows that the composition
\[
    \delta = \delta_1^{-1} \circ\left({\frac{\partial~~}{\partial t^k} } \right)\circ \delta_2\,.
\]
is a Legendre field between the two generalised Frobenius manifolds constructed from $FM_{1,2}$ via the Legendre fields $\delta_{1,2}\,.$ To find its homogeneity we first introduce the notation $\mathcal{L}_E\delta_i = [\delta_i]\delta_i$ for the degrees $[\delta_i]$ of the homogeneous Legendre fields $\delta_i$ and note that ${\frac{\partial~~}{\partial t^k} }$ has degree $-d_k\,.$
Proposition \ref{prop:spectrumhomogeneousLF}(i) gives the following relations between the various charges:
\begin{eqnarray*}
    \tilde{d} & = & d - 2 - 2 [\delta_1]\,,\\
    {\rm charge~of~}{\rm FM}_2 & = & d - 2 + 2 d_k\,.\\
    \widehat{d} & = & {\rm charge~of~}{\rm FM}_2 -2 - 2 [\delta_2]\,,\\
    \widehat{d} & = & \tilde{d} - 2 - 2[\delta]\,.
\end{eqnarray*}
and from its definition, 
\begin{eqnarray*}
[\delta] & = & [\delta_1^{-1}] + (-d_k) + [ \delta_2] +2 \,,\\
& = & -d_k - [\delta_1] + [\delta_2]\,.
\end{eqnarray*}
It follows trivially that if $\tilde{d}=\widehat{d}$ then $[\delta]=-1$ and if this holds then the degrees of $\delta_1$ and $\delta_2$ have to be related by the condition $[\delta_2]=[\delta_1] - 1 + d_k\,.$

\item[(b)] If the connecting Legendre transformation exists, then since Legendre transformations preserve the spectrum (by Corollary \ref{Cor31}), the two generalized Frobenius manifolds must have the same spectrum. Conversely, if ${\rm gFM_1}$ and ${\rm gFM_2}$ have the same spectrum then the charges of the corresponding Frobenius manifolds are given directly (again by Corollary \ref{Cor31}) in terms of elements of the spectrum, namely
\begin{eqnarray*}
    d_{{\rm FM}_1} & = & - 2 \mu_\bullet\,,\\
    d_{{\rm FM}_2} & = & - 2 \mu_\circ\,,
\end{eqnarray*}
for some $\mu_\bullet\,,\mu_\circ$ in the (common) spectrum. The degree of any connecting Legendre field, if it exists, must then satisfy
\[
(- 2 \mu_\circ) = (- 2 \mu_\bullet) - 2 - 2 [\delta]
\]
or $[\delta]=\mu_\circ-\mu_\bullet -1$ and it is easy to see that a unique flat vector field $\partial_{t^i}$ for ${\rm FM}_1$ must exist with this homogeneity - namely the one with $\mu_i = \mu_\circ$ which exists as the spectrum is common to all the manifolds. 
 
\end{itemize}

\end{proof}

\noindent Note that, if $k=n$ then from the symmetries of the spectrum, $d_n=1-d$ and hence the map in part (a) becomes
\[
    \mathbb{P}_{d,\tilde{d}}^{n-1} \xlongrightarrow{\hspace{3mm}\frac{\partial~~}{\partial t^n} \hspace{3mm}} \mathbb{P}_{-d,\tilde{d}}^{n-1}
\]
(this will be shown explicitly in the case $n=2$ in Section 4). More generally, there will always exist pairs with $\mu_\bullet+\mu_\circ =0 $ which will give Frobenius manifolds ${\rm FM}_{1,2}$ with opposite charge. 
Finally, in part (b), if ${\rm FM}_1\cong {\rm FM_2}$ then the connecting Legendre field exists - it must be the identity field and hence $[\delta]=-1$ which implies that $\mu_\bullet=\mu_\circ$ and hence elements of the spectrum are not all distinct. In two-dimensions this implies, since $\mu_\bullet+\mu_\circ =0\,,$ that $\mu_\circ=\mu_\bullet$ and hence that these Frobenius manifolds have charge zero.

\section{Classification of 2D homogeneous Legendre fields}\label{classification}


In this section we classify all homogeneous semi-simple Legendre fields in 2-dimensions by reducing - on using homogeneity - the defining equations for the Legendre field to second order ordinary differential equations for the potentials, and these equations turn out to be the hypergeometric equations in all cases. These calculations may be done in both flat coordinates and canonical coordinates.

Following \cite{Dubrovin1996}, there is an equivalence relation on the set of Frobenius manifolds. Two Frobenius manifolds are equivalent if there is an affine-linear map in the flat coordinates (equivalently, a diffeomorphism) which is an isomorphism of the algebras in the tangent spaces, and which rescales the metric by a non-zero constant. 
Equivalence classes of two-dimensional semi-simple Frobenius manifolds with isotropic identity are parametrised by the charge which can take arbitrary complex values, and we will identify such equivalent Frobenius manifolds. The resulting list of Frobenius manifolds is given in Table 1.

Legendre transformations of two-dimensional semisimple Frobenius manifolds with non-zero charge do not generate such equivalences, excluding the trivial case of the Legendre field being a scalar multiple of the identity field $e$. Indeed, if this were the case, one would have $\widetilde{\nabla}=\nabla$ so a Legendre field $\delta$ would have to satisfy the equation $\delta\circ\nabla_XY=\nabla_X(\delta\circ Y)$ and therefore it would have to be flat. Since $\delta$ is the unit of the new Legendre transformed Frobenius manifold it has to be isotropic, and the charge should be preserved. One can check that in the case of non-zero charge this implies that $\delta$ is proportional to the identity $e$ rather than to the other flat coordinate direction.

\subsection{Flat coordinates}

To fix notation, the prepotentials, Euler field and the charge of these Frobenius manifolds are given in Table 1 (we exclude the single case when $\eta(e,e)=1$ which yields a cubic prepotential. For details, see \cite{LeoThesis}). While $FM_1$ is defined for any value of $N$, if
$N \in\mathbb{N}_{\geq 3}$ then the corresponding Frobenius manifold has the structure of a dihedral group orbit space $\mathbb{C}^2/I_2(N)\,$ and such solutions have particularly nice algebraic properties. For other values of $N$ one has multi-valued functions with branching.

\begin{table}
  \[
  \begin{array}{llll}
  {FM_1:~} & \displaystyle{F=\frac{1}{2} (t^1)^2 t^2 + \frac{(t^2)^{N+1}}{N(N^2-1)}}\,,& \displaystyle{E=t^1\frac{\partial~~}{\partial t^1} + \frac{2}{N} t^2\frac{\partial~~}{\partial t^2}}\,, & {\displaystyle d= 1- \frac{2}{N} }\,,\,(N \neq 0,\pm 1)\,,\\
  
&&&\\
{FM_2:~} & \displaystyle{F=\frac{1}{2} (t^1)^2 t^2 + \frac{1}{2} (t^2)^2 \log t^2}\,,& \displaystyle{E=t^1\frac{\partial~~}{\partial t^1} + 2t^2\frac{\partial~~}{\partial t^2}}\,, & {\displaystyle d= -1 }\,,\\
&&&\\
 
{FM_3:~} & \displaystyle{F=\frac{1}{2} (t^1)^2 t^2 + \frac{1}{2} \log t^2}\,,& \displaystyle{E=t^1\frac{\partial~~}{\partial t^1} - 2t^2\frac{\partial~~}{\partial t^2}}\,, & {\displaystyle d= 3 }\,,\\

&&&\\

{FM_4:~} & \displaystyle{F=\frac{1}{2} (t^1)^2 t^2 + e^{t^2} }\,,& \displaystyle{E=t^1\frac{\partial~~}{\partial t^1} + 2\frac{\partial~~}{\partial t^2}}\,, & {\displaystyle d= 1 }\,.\\  
&&&\\

  \end{array}
  \]
  \caption{The semi-simple Frobenius manifolds in 2-dimensions with antidiagonal metric (equivalently $\eta(e,e)=0).$}
\end{table}

We first rewrite the homogeneity condition on the potential $h$ in terms of the pair $\{d,\tilde{d}\}\,,$

\[
\mathcal{L}_E h = 
    \begin{cases}
        \left(1- \frac{1}{2} (d+\tilde{d}) \right) \, h\,, & {\rm if~} d+\tilde{d} \neq 2\,,\\
       {c} \,, &{\rm if~} d+\tilde{d}  = 2\\
    \end{cases}\]
where $c$ is a constant, and from the Euler vector fields given on Table 1 one can construct a scaling ansatz for the potential. Since all these manifolds have an anti-diagonal metric, the Legendre field equation reduces to the partial differential equation
\begin{equation}\label{LFeqn}
\frac{\partial^2 h}{\partial (t^2)^2} = f^{\prime\prime\prime}(t^2) \,\frac{\partial^2 h}{\partial (t^1)^2}
\end{equation}
where $f(t^2)$ is the term independent of $t^1$ in the prepotential of the Frobenius manifold, and with the scaling anstaz this then reduces to a second order linear equation. The analysis splits depending on the value of $d+\tilde{d}\,.$

\begin{table}
  \[
  \begin{array}{llll}
  {FM_1:~} 
    & \displaystyle{h=(t^2)^\alpha \,\,y\left( \frac{N^2 (t^1)^2}{4 (t^2)^N}\right)}\,,&&\alpha=\frac{1}{4}(2+N-N\tilde{d})\,;\\
&&&\\
{FM_2:~} 
 & \displaystyle{h=(t^2)^\alpha \,\,y\left( \frac{(t^1)^2}{4 t^2}\right)}\,,&& \alpha=\frac{1}{4}(3-\tilde{d})\,;\\
    &\\
{FM_3:~} & \displaystyle{h=(t^2)^\alpha \,\,y\left( \frac{(t^1)^2 t^2}{4}\right)}\,,&& \alpha = \frac{1}{4}(1+\tilde{d})\,;\\
    &\\
{FM_4:~} 
& \displaystyle{h=e^{\alpha t^2} \,\,y\left( \frac{(t^1)^2}{4 e^{t^2}}\right)}\,,&& \alpha = \frac{1}{4}(1-\tilde{d})\,.\\
    &\\
&&&\\
  \end{array}
  \]
  \caption{The scaling ansatz for the Legendre field potential in term of the charge $\tilde{d}\,.$}
\end{table}

\subsubsection{The case $d+\tilde{d}\neq 2\,.$}

The scaling ansatz for the potentials are given in Table 2. With these, equation (\ref{LFeqn}) reduces to the hypergeometric equation
\[
x(1-x) y^{\prime\prime}(x) + (c-(a+b+1)x) \,y^{\prime}(x) - a \, b\, y(x) = 0\
\]
with the values of the constants $a\,,b$ and $c$ being given in Table 3, this also containing the values of the constants $\lambda\,,\mu$ and $\nu$ for the corresponding Schwarz triangles (whose angles are $\lambda\pi\,,\mu\pi$ and $\nu\pi$). Note that the singular points of the hypergeometric equation correspond to interesting strata within the Frobenius manifold.

Given the simple relationship between $\alpha$ and $\tilde{d}$ (see Table 2), we get for each $\tilde{d}$ a solution which is a linear combination of two independent solutions to the hypergeometric equation. Since the Legendre field must be invertible (and hence non-zero), and an overall scale plays no role, we obtain the solution space $\mathbb{P}^1_{d,\tilde{d}}$ (isomorphic to $\mathbb{CP}^1$ for all values of $d$ and $\tilde{d}$) of generalized Frobenius manifolds that can be constructed from a given 2D-Frobenius manifold by this construction.

Considerable care has to be taken in constructing global solutions of these hypergeometric equations and hence global Legendre fields.  While the elementary solution (where $\alpha_1$ and $\alpha_2$ are constants)
\[
y(x) = \alpha_1 \,\,\, {}_2F_1(a,b,c;x) + \alpha_2 \,\,\, \,\,x^{1-c}\,{}_2F_1(a-c+1,b-c+1,2-c;x)
\]
exists in the region $|x|<1$ around the point $x=0$, the parameter values are such that the standard connection formulae may not hold (for example, for $FM_{2,3,4}$ we have $a-b=0$ or $1$ and standard formulae connecting solutions around $x=0$ to $x=\infty$ only hold for $a-b\notin\mathbb{Z}$). Instead one must use more subtle connection formulae \cite{haraoka} and full details, including explicit formulae for the Legendre fields for all values of $\alpha$ may be found in \cite{LeoThesis}.

\begin{center}
\begin{table}
\vskip 2mm
\begin{tabular}{c ||c |c| c || c |c |c} 
& $a$ & $b$ & $c$ & $\lambda=|1-c|$ & $\nu=|b-a|$ & $\mu=|c-a-b|$ \\
\hline
&&&&&&\\
{$FM_1$} & $\displaystyle{-\frac{\alpha}{N}}$ & $\displaystyle{\frac{1}{N}-\frac{\alpha}{N}}$ & $\displaystyle{\frac{1}{2}}$ &$\displaystyle{\frac{1}{2}}$&$\displaystyle{\frac{1}{N}}$&$\displaystyle{\left|\frac{1}{2} + \frac{2\alpha-1}{N}\right|}$\\
 &&&&&&\\
 {$FM_2$} & $\displaystyle{-\alpha}$ & $\displaystyle{1-\alpha}$ & $\displaystyle{\frac{1}{2}}$ &$\displaystyle{\frac{1}{2}}$&$1$&$\displaystyle{\left|\frac{1}{2} - 2\alpha\right|}$\\
 &&&&&&\\
 {$FM_3$} & $\displaystyle{\alpha}$ & $\displaystyle{\alpha-1}$ & $\displaystyle{\frac{1}{2}}$ &$\displaystyle{\frac{1}{2}}$&$1$&$\displaystyle{\left|\frac{3}{2} - 2 \alpha\right|}$\\
 &&&&&&\\
 {$FM_4$} & $\displaystyle{-\alpha}$ & $\displaystyle{-\alpha}$ & $\displaystyle{\frac{1}{2}}$ &$\displaystyle{\frac{1}{2}}$&$0$&$\displaystyle{\left|\frac{1}{2} + 2\alpha\right|}$\\
\end{tabular}
\vskip 5mm
\caption{Values of the hypergeometric equation parameters and the Schwarz triangle angles $\lambda\pi\,,\mu\pi\,,\nu\pi\,.$}
\end{table}
\end{center}

Given the transcendental nature of these solutions, one could impose a requirement of a \lq suitable algebraic condition\rq~on the solution, such as the solution being an algebraic function. This implies that the solution has finite monodromy, and as the monodromy group
for the hypergeometric equation depends on the value of the constant $\alpha$ one just has to compare the values of $\lambda\,,\mu$ and $\nu$ against those in Schwarz's list of algebraic solutions of the hypergeometric equation and this then isolates the corresponding values of $\alpha\,.$ One can also extend this analysis to obtain the conditions on $\alpha$ required for these angles to define a triangle groups. As expected, if $\tilde{d}=1$ then one obtains interesting finite monodromy and, more generally, triangle groups, namely $\Delta(2,2,N), \,\Delta(1,2,2)$ (twice) and $\Delta(2,2,\infty)\,.$ 
However, in each case there are other specific values which, from Schwarz's list, produces a finite monodromy group. 

\begin{ex}\label{FM1example}
    \noindent{$FM_1$:} Here we restrict $N$ to be a positive integer greater or equal to three. For these values the manifold itself is orbit space $\mathbb{C}^2/I_2(N)\,.$ From Schwarz's list one can then find values of the constant $\alpha$ for which the hypergeometric equation has finite monodromy. This isolates, for arbitrary $N$, the value $\alpha=\tfrac12$ with the corresponding angles $\left\{ \tfrac{\pi}{2},\tfrac{\pi}{2},\tfrac{\pi}{N} \right\}$ and dihedral monodromy group $\Delta(2,2,N).$ Other values of $\alpha\,,$ restricting to $0<\alpha\leq \frac{1}{2}\,,$ for which the monodromy group is finite are given by:
\[
\begin{array}{lll}
N=3\,, & \alpha=\tfrac14\,, & \bigl\{\tfrac12,\tfrac13,\tfrac13\bigr\}\,;\\[0.5em]
&\alpha=\tfrac18\,, & \bigl\{\tfrac12,\tfrac13,\tfrac14\bigr\}\,;\\[0.5em]
&\alpha=\tfrac{1}{20}\,, & \bigl\{\tfrac12,\tfrac13,\tfrac15\bigr\}\,;\\[0.5em]
&\alpha=\tfrac{7}{20}\,, & \bigl\{\tfrac12,\tfrac13,\tfrac25\bigr\}\,;\\[0.5em]
N=4\,,& \alpha = \tfrac16\,, & \bigl\{\tfrac12,\tfrac13,\tfrac14\bigr\}\,;\\[0.5em]
N=5\,,& \alpha = \tfrac{1}{4}\,, & \bigl\{\tfrac12,\tfrac15,\tfrac25\bigr\}\,;\\[0.5em]
& \alpha = \tfrac{1}{12}\,, & \bigl\{\tfrac12,\tfrac13,\tfrac15\bigr\}\,.
\end{array}
\]
\end{ex}

\medskip

\begin{ex}
\noindent{$FM_{2,3,4}$:} The remaining cases may be analyzed in a similar fashion and we omit the full details. Using the fact that exponent differences are only considered modulo integers, the cases $FM_{2,3,4}$ give, if $\alpha = (2-m)/(4m)$, where $m$ is an integer, the triangle group $\Delta(2,m,\infty)\,.$ If $\alpha=0$ (see below) then this is the Euclidean triangle group $\Delta(2,2,\infty)$ rather than (for $m\geq 3$) a hyperbolic triangle group.
\end{ex}

\subsubsection{The case $d+\tilde{d}= 2\,.$}

In this case the scaling ansatz, for each two-dimensional Frobenius manifolds, is just given via the replacing 
\[
h \longmapsto c \log(t^1) + \left.h\right|_{\alpha=0}
\]
for each $h$ in Table 2. On substituting these into equation (\ref{LFeqn}) one obtains, in each case, an inhomogeneous hypergeometric equation with the inhomogeneous term depending on the arbitrary constant $c$. The homogeneous equations are the same hypergeometric equations as above but with the parameter values corresponding to the condition $\alpha=0$. This implies (from Table 3) that in all cases $a=0$ and hence these equations are just first-order linear equations in $y'(x)$ and hence these can be solved by quadrature. We omit the details, and the final results may be found in \cite{LeoThesis}. From Table 2 it also follows that if $\alpha=0$ then one will only obtain a canonical generalized Frobenius manifold (which, by definition, has $\tilde{d}=1$) for the Frobenius manifold $FM_4\,.$

\begin{ex}
    \noindent{$FM_4$:} With
\[
    h= c \log( t^1) + y\left( \frac{(t^1)^2}{4 e^{t^2}}\right)
    \]
one obtains the differential equation
    \[
    z(1-z) y^{\prime\prime}(x) + \left(\frac{1}{2} - z\right) y^\prime(x) = \frac{c}{4z}
    \]
which may easily be solved by quadrature. The solution depends on three constants but one just enters as an additive constant and may be ignored as it does not contribute to the Legendre field. From $\delta = (dh)^\sharp$ this only gives a two-parameter family of solutions in agreement with Example \ref{deltaP1}\,. 
    \end{ex}

\medskip

\subsection{Canonical coordinates}

Two-dimensional semi-simple Frobenius manifolds are characterized solely by their charge $d\,;$ the various different forms in Table 1 result from the integrations required to find the flat coordinates and prepotential, this generating a number of special cases. In this section we will denote by $FM(d)$ the 2D-Frobenius manifold with charge $d\,.$

If a Frobenius manifold is semi-simple, the Darboux-Egorov equations for the rotation coefficients remain unchanged under a Legendre transformation - both the \FM and the \gFM have the same rotation coefficients. What does change is how the metric is reconstructed from these coefficients \cite{SS}. 

Recall, on writing $\eta=\sum H_i^2(u)\, du_i^2$ that the rotation coefficients are
\[
\beta_{ij} = \frac{\partial_i H_j}{H_i}
\]
and that for 2D-Frobenius manifolds,
\[
\beta_{12}=\beta_{21} = \frac{-i \mu}{u_1-u_2}\,, \qquad \mu=-\frac{d}{2}\,,
\]
and that the requirements on the metric imply $\mathcal{L}_E(H_i)=\mu H_i$ and $\mathcal{L}_e H_i =0$ and hence
\begin{eqnarray*}
    H_1 & = & \phantom{-i}(u_1-u_2)^\mu\,, \\
    H_2 & = & -i(u_1-u_2)^\mu\,.
\end{eqnarray*}
For generalized Frobenius manifolds one drops the last requirement which, in 2-dimensions, imply that $H_i=H_i(u_1-u_2)\,.$ Thus on defining
\[
\widetilde{\eta} = \sum_i {\widetilde{H}}_i^2 \, du_i^2\qquad ({\rm so~} \delta = \sum_i \frac{\widetilde{H}_i}{H_i} \frac{\partial~~}{\partial u_i})
\]
one obtains the system
\begin{eqnarray*}
    \frac{\partial \widetilde{H}_1}{\partial u_2 } & = & 
    \frac{-i\mu}{u_1-u_2} {\widetilde{H}}_2\,,\\
    \frac{\partial \widetilde{H}_2}{\partial u_1 } & = & \frac{-i\mu}{u_1-u_2} {\widetilde{H}}_1\,,
    \end{eqnarray*}
with $\mu=-d/2$ and the single constraint $\mathcal{L}_E({\widetilde{H}_i)=-\displaystyle{\frac{\tilde{d}}{2}} \,{\widetilde{H}}_i}\,.$

On eliminating ${\widetilde{H}}_2$ one obtains a Euler-Poisson-Darboux type equation
\[
(u_1-u_2) \frac{\partial~}{\partial u_1} \left( 
(u_1-u_2)\frac{\partial\widetilde{H}_1}{\partial u_2}
\right) + \mu^2\,\widetilde{H}_1 = 0
\]
and with the scaling ansatz
\[
\frac{{\widetilde{H}}_1}{H_1} = u_1^{-\frac{{\tilde{d}}}{2}} u_2^{+\frac{d}{2}}f\left(\frac{u_1}{u_2}\right)
\]
one obtains a hypergeometric equation for $f(x)$ with parameter values $a=-d/2\,, b=1-d/2-\tilde{d}/2\,,c=1-\tilde{d}/2\,$ whose monodromy group depends $d$ and $\tilde{d}\,.$  Note that the singular points of the hypergeometric equation correspond to points on the discriminant $u_1=0$ or $u_2=0$ (the regular singular points $x=0$ and $x=\infty$ respectively) or on the caustic $u_1=u_2$ (the regular singular point $x=1$).

The solution to this hypergeometric equation may be written using Riemann's $P$-notation as
\[
f(x) = P \left\{
\begin{array}{ccc|c}
0 & 1 & \infty & \\
0 & 0 & -d/2 & x \\
\tilde{d}/{2} & d & 1-d/2-\tilde{d}/2 & 
\end{array} 
\right\}
\]
\noindent and around the point $z=0$ the solution may be written (for generic values of $d$ and $\tilde{d}$) in terms of hypergeometric functions, 
\[
f(x) = c_1\,\, {}_2F_1\left(-\frac{d}{2}, 1-\frac{d}{2}-\frac{\tilde{d}}{2}, 1-\frac{\tilde{d}}{2};x\right) + c_2 \,\, x^{\tilde{d}/2} \, {}_2F_1\left(1-\frac{d}{2},-\frac{d}{2}+\frac{\tilde{d}}{2},1+\frac{\tilde{d}}{2},x\right)
\]
where $c_1$ and $c_2$ are constants (for non-generic values we refer the reader to the standard references\footnote{For example: NIST Digital Library of Mathematical Functions, https://dlmf.nist.gov .}). To obtain $\tilde{H}_2$ one just uses the first equation of the above system (similar results also appear in \cite{AL13}). To find the Legendre field in flat coordinates one just has to perform the appropriate changes of variables. Since canonical coordinates, in two dimensions, are obtained by solving the quadratic equation $\det(g^{-1} - u \eta^{-1}) = 0$ this is a straightforward, but tedious, exercise. An alternative, but equivalent approach, is to directly solve the Legendre field equations in flat-coordinates for each of the 2-dimensional Frobenius manifolds, see \cite{LeoThesis}. Further results on properties of canonical coordinate description of Legendre fields may be found in \cite{JohanThesis}.

Recall, for two-dimensional semi-simple Frobenius manifolds, that Theorem \ref{mapbetweenspaces} gives a map between the solution spaces $\mathbb{P}^1_{d,\tilde{d}}$ and $\mathbb{P}^1_{-d,\tilde{d}}\,.$ Since these spaces are the solution spaces to hypergeometric equations, this map should correspond to some property of the hypergeometric equation itself. Recall the following property of solutions of the hypergeometric equation:
\[
P \left\{
\begin{array}{ccc|c}
0 & 1 & \infty & \\
0 & 0 & a & z \\
1-c & c-a-b & b & 
\end{array} 
\right\}= (1-z)^{c-a-b} \,\,
P \left\{
\begin{array}{ccc|c}
0 & 1 & \infty & \\
0 & 0 & c-a & z \\
1-c & a+b-c & c-b & 
\end{array} 
\right\}\,.
\]
Rewriting this with the above values of $a\,,b$ and $c$ gives
\[
P \left\{
\begin{array}{ccc|c}
0 & 1 & \infty & \\
0 & 0 & -d/2 & x \\
\tilde{d}/{2} & d & 1-d/2-\tilde{d}/2 & 
\end{array} 
\right\} = (1-z)^d \,\, P \left\{
\begin{array}{ccc|c}
0 & 1 & \infty & \\
0 & 0 & d/2 & x \\
\tilde{d}/{2} & -d & 1+d/2-\tilde{d}/2 & 
\end{array} 
\right\}
\]
which shows how the solutions with $\pm d$ are related. The term $(u_1-u_2)^d$ comes from expressing the connecting Legendre field in terms of canonical coordinates, i.e.
\[
\frac{\partial~~}{\partial t^2} = 2 (u_1-u_2)^d \left(\frac{\partial~~}{\partial u^1} - \frac{\partial~~}{\partial u^2}\right)\,.
\]

Canonical generalized Frobenius manifolds have charge $\tilde{d}=1$ and for this value the hypergeometric equation has a closed-form solution for all values of $d\,,$ namely
\[
f(x) = c_1 (1+\sqrt{x})^d + c_2 (1-\sqrt{x})^d\,.
\]
with monodromy group - generically an infinite group - depending solely on the value of $d\,.$ Reconstructing the Legendre field $\delta = \delta^1(u)\, \partial_{u_1} + \delta^2(u)\, \partial_{u_2}$ gives
\begin{equation}\label{deltad}
\begin{array}{ccc}
    \delta^1(u) & = & \frac{1}{\sqrt{u_1}} \left\{
    c_1 (\sqrt{u_1}+\sqrt{u_2})^d +c_2 (\sqrt{u_1}-\sqrt{u_2})^d\right\})\,,\\
    &&\\
    \delta^2(u) & = & \frac{1}{\sqrt{u_2}} \left\{
    c_1 (\sqrt{u_1}+\sqrt{u_2})^d -c_2 (\sqrt{u_1}-\sqrt{u_2})^d\right\})
    \end{array}
\end{equation}
which again shows the branching along the discriminant. 

Drawing these results together gives:

\begin{thm} 
Let $FM(d)$ be a two-dimensional semi-simple Frobenius manifold of charge $d$ and let $\delta$ be a homogeneous Legendre field for this Frobenius manifold, with $\mathcal{L}_E \delta = \kappa \delta$. In canonical coordinates, the components of $\delta$ are given in terms of solutions of the hypergeometric equation with parameters
\[
a=-d/2\,, \qquad b=1-d/2-\tilde{d}/2\,,\qquad c=1-\tilde{d}/2\,
\]
where $\tilde{d}=d-2 - 2 \kappa\,.$

The map from the 2-dimensional semi-simple Frobenius manifold $FM(d)$ to the corresponding canonical generalized Frobenius manifold (so $\tilde{d}=1$) is given by the Legendre field
\[
\delta_d = \delta^1(u)\, \partial_{u_1} + \delta^2(u)\, \partial_{u_2}\]
where 
\begin{eqnarray*}
    \delta^1(u) & = & \frac{1}{\sqrt{u_1}} \left\{
    c_1 (\sqrt{u_1}+\sqrt{u_2})^d +c_2 (\sqrt{u_1}-\sqrt{u_2})^d\right\})\,,\\
    \delta^2(u) & = & \frac{1}{\sqrt{u_2}} \left\{
    c_1 (\sqrt{u_1}+\sqrt{u_2})^d -c_2 (\sqrt{u_1}-\sqrt{u_2})^d\right\})\,
\end{eqnarray*}
with $c_1$ and $c_2$ constants.
\end{thm}

\noindent Finally, given a 2-dimensional generalized Frobenius manifold one can construct, by performing Legendre transforms using the two flat homogeneous Legendre fields (the two coordinate vector fields), two different Frobenius manifolds. It also follows from the transitive action of composition of Legendre field that these two Frobenius manifolds are related by a flat-Legendre field, and this construction changes of the sign of the charge.

\subsection{The degenerate case $d=0$}

Throughout this section there has been a tacit assumption that $d\neq 0\,.$ If $d=0$ then the prepotential and Euler field takes the form (see Table 1)
\[
F = \frac{1}{2} (t^1)^2 t^2 + \frac{1}{6} (t^2)^3\,, \qquad E = t^1 \frac{\partial~~}{\partial t^1 } + t^2 \frac{\partial~~}{\partial t^2 } 
\]
(and conversely, if $d_1=d_2\,,$ then the standard normalization requires $d_1=1$ and equation (\ref{muconstraint}) gives $d=0$). A simple linear change of variable shows that this is equivalent to the direct sum of two one-dimensional Frobenius manifolds
\[
\tilde{F} = \frac{1}{6} (\tilde{t}^1)^3 + \frac{1}{6} (\tilde{t}^2)^3\,, \qquad E = {\tilde t}^1 \frac{\partial~~}{\partial {\tilde t}^1 } + {\tilde t}^2 \frac{\partial~~}{\partial {\tilde t}^2 } 
\]
with metric $\tilde{\eta} = (d\tilde{t^1})^2 - (d\tilde{t^2})^2\,$ (these  coordinates also serve as canonical coordinates). The Legendre field equations (\ref{defLF}) then decouple to give, with the assumption of homogeneity, the Legendre field
\[
\delta = c_1 \,({\tilde t}^1)^{\kappa+1} \frac{\partial~~}{\partial {\tilde t}^1 } + c_2\, ({\tilde t}^2)^{\kappa+1} \frac{\partial~~}{\partial {\tilde t}^2 } 
\]
where $c_1$ and $c_2$ are constants. To obtain a canonical generalized Frobenius manifold requires $\kappa = -3/2$ and this agrees with the explicit formulae given above with $d=0\,.$

\section{Examples}

In this section we examine the known examples of 2-dimensional generalized Frobenius manifolds can see how they fit into the theory developed above. Since applying a flat Legendre field to a \gFM results in a \FM we first apply this to obtain the underlying Frobenius manifold. The construction involves inverting a change of variable to obtain $\tilde{t}^i=\tilde{t}^i({\bf t})$ and since $\tilde{t}^n=h$
this will automatically give the potential for the (non-flat) Legendre field that takes one from the \FM to the generalized Frobenius manifold with which one started. This will enable the specific point in the space $\mathbb{P}^1_{d,\tilde{d}}$ to be identified.

\subsection{Examples reducing to Coxeter orbit space $\mathbb{C}^2/I_2(3)$}\label{example51}

\begin{ex}
Starting with the solution found in \cite{orbit1} (suitably scaled) on the affine orbit space $\mathbb{C}^2/{A_2^{aff}}\,:$
\begin{eqnarray*}
\widetilde{F} & = & \frac{1}{240} {(\tilde{t^2})^6} - 
\frac{i}{24} (\tilde{t}^2)^4 \tilde{t}^1 -
\frac{1}{4} (\tilde{t}^2 \tilde{t}^1)^2 + 
\frac{i}{6} {(\tilde{t^1})^3}\,,\\
&\\
E & = & \frac{2}{3} \tilde{t}^1 \frac{\partial~~}{\partial {\tilde{t}^1}} + \frac{1}{3} \tilde{t}^2 \frac{\partial~~}{\partial {\tilde{t}^2}}
\end{eqnarray*}
(and $\tilde{d}=1$) and applying the flat Legendre field
\[
\delta^{-1} = \frac{\partial~~}{\partial {\tilde{t}^1}}
\]
results in the change of variables
\begin{eqnarray*}
    t^1 & = & {\widetilde{F}_{12}} = - \tilde{t}^1 \tilde{t}^2 - \frac{i}{6} (\widetilde{t}^2)^3\,,\\
    t^2 & = & {\tilde{F}_{11}} = i\, \tilde{t}^1 - \frac{1}{2} (\tilde{t}^2)^2
\end{eqnarray*}
and the equation
\begin{eqnarray*}
    \widetilde{F}_{22} & = & \frac{1}{2} \left\{ 
    i \tilde{t}^1 -  \frac{1}{2} (\tilde{t}^1)^2\right\}^2 \,,\\
    & = & \frac{1}{2} (\widetilde{F}_{11})^2
\end{eqnarray*}
and hence the resulting \FM is the orbit space $\mathbb{C}^2/I_2(3)\,.$ Eliminating $\tilde{t^1}$ give the cubic
\[
\frac{1}{3}
 (\tilde{t}^2)^3 + t^2 \tilde{t}^2 + i t^1 = 0\,.
\]
This may be solved using Cardano's formula to give (since $h=\tilde{t}^2\,,$ and up to an overall constant)
\[
h_1(t^1,t^2) = \left(
t^1 + \sqrt{(t^1)^2 - \frac{4 (t^2)^3}{9} }
\right)^\frac{1}{3} + \left(
t^1 - \sqrt{(t^1)^2 - \frac{4 (t^2)^3}{9} }
\right)^\frac{1}{3}
\]
which corresponds to the point $\alpha_1=\alpha_2\,$ (or $[1,1] \in \mathbb{P}^1_{d=\frac{1}{3},\tilde{d}=1}$ in projective coordinates) in the general solution found above.
\end{ex}

\begin{ex}

Starting with the solution found in \cite{orbit1} (suitably scaled) on the affine orbit space $\mathbb{C}^2/{G_2^{aff}}\,$:
\begin{eqnarray*}
\widetilde{F} & = & \frac{1}{960} {(\tilde{t^2})^6} - 
\frac{1}{8} (\tilde{t}^2 \tilde{t}^1)^2 + 
\frac{(\tilde{t^1})^4}{12(\tilde{t}^2)^2}\,,\\
&\\
E & = & \frac{2}{3} \tilde{t}^1 \frac{\partial~~}{\partial {\tilde{t}^1}} + \frac{1}{3} \tilde{t}^2 \frac{\partial~~}{\partial {\tilde{t}^2}}
\end{eqnarray*}
(and $\tilde{d}=1$) and applying the flat Legendre field
\[
\delta^{-1} = \frac{\partial~~}{\partial {\tilde{t}^1}}
\]
results after a similar calculation to the potential
\[
h_2(t^1,t^2) = \left(
t^1 + \sqrt{(t^1)^2 - \frac{4 (t^2)^3}{9} }
\right)^\frac{1}{3} - \left(
t^1 - \sqrt{(t^1)^2 - \frac{4 (t^2)^3}{9} }
\right)^\frac{1}{3}
\]
which gives the condition $\alpha_1=-\alpha_2$ (or $[1,-1] \in \mathbb{P}^1_{\tilde{d}=1}$ in projective coordinates).
\end{ex}

The functions $h_1$ and $h_2$ are the potentials for the Legendre fields $\delta_1$ and $\delta_2$ in Figure \ref{fig:fig4}.
A third solution, corresponding to the case $\alpha_1=0$ (or $\alpha_2=0$) (or $[1,0] {\rm ~or~} [0,1] \in \mathbb{P}^1_{d=\frac{1}{3},\tilde{d}=1}$ in projective coordinates) will be given in the next section.

\subsection{Examples reducing to extended affine Weyl orbit space $\mathbb{C}^2/A^{(1)}_1$}

That fact that applying a (flat)-Legendre transformation to the generalized Frobenius manifolds defined by
\begin{equation}
\begin{array}{ccl}
    \widetilde{F} & = & \displaystyle{\frac{1}{2}
(\tilde{t}^1)^2 \tilde{t}^2 + \tilde{t}^1 e^{\tilde{t}^2} + \frac{1}{2} (\tilde{t}^1)^2 \log \tilde{t}^1}\,,\\
&&\\
E & = & \displaystyle{\tilde{t}^1 \frac{\partial~~}{\partial \tilde{t}^1} + 
\ \frac{\partial~~}{\partial \tilde{t}^2}}
\end{array}\label{gQP1}
\end{equation}
(with $\delta^{-1}=\frac{\partial~~}{\partial \tilde{t}^1}$) results in the Frobenius manifold structure on the extended affine Weyl orbit space $\mathbb{C}^2/A^{(1)}_1$
has been known for some time. It formed the basis of the connection between the quantum cohomology of certain resolved conifolds \cite{antidiag} and played a crucial role in the development of the concept of a generalized Frobenius manifold. Applying the Legendre transformation $\delta^{-1}=\frac{\partial~~}{\partial \tilde{t}^2}$ results in the Frobenius manifold $FM_4$, reflecting that fact that $FM_2$ and $FM_4$ are themselves related by a flat-Legendre field.

The general form of the Legendre field $\delta$ {\sl from} the Frobenius manifold $FM_4$ was derived in Example \ref{deltaP1} above.
By coordinate transformation similar to the previous example one may find the corresponding point in the space $\mathbb{P}^1_{d=1,\tilde{d}=1}$:
\begin{eqnarray*}
[1,0] \in \mathbb{P}^1_{d=1,\tilde{d}=1} &\Rightarrow& FM_4 ~~~{\rm (identity)}\\
{}[1,1] \in \mathbb{P}^1_{d=1,\tilde{d}=1} &\Rightarrow& (\ref{gQP1})
\end{eqnarray*}
By direct calculation one may show that at the point $[0,1] \in \mathbb{P}^1_{d,\tilde{d}=1}$ in the solution space one can invert the change of variable, and this results in the new generalized Frobenius manifold
\begin{equation}
    \widetilde{F}  =  \frac{1}{2}
(\tilde{t}^1)^2 \tilde{t}^2  +  (\tilde{t}^1)^2 \log \tilde{t}^1 - (\tilde{t}_1)^2 \log[ 1 + e^{\tilde{t}^2}]\,,
\label{QP1-01}
\end{equation}
with
\[
\tilde{e}^\flat = d \left[ \log \left(
\frac{ e^{\tilde{t}^2}}{\tilde{t^1}-e^{\tilde{t}^2}}\right)\right]\,.
\]
With arbitrary $[\alpha_1,\alpha_2]$ one obtains, geometrically (but not algebraically as one cannot invert the change of variables), an interpolating family of generalized Frobenius manifolds connecting $FM_4$ to $(\ref{QP1-01})$ via (\ref{gQP1}).

\section{The solution space of 2D generalized Frobenius manifolds}

For Frobenius manifolds, the extra condition $\nabla e=0$ ties together the prepotential and the metric in a very strong way. In flat coordinates one takes $e=\partial_1$ and the normalization $\eta_{ij}=c_{1ij}$ determines the $t_1$-dependence of the prepotential. This generalises to  situations the metric may be defined as an arbitrary linear combination of third-derivatives of the prepotential (see, for example \cite{Martini}). For generalized Frobenius manifold the link between the prepotential and the metric is weakened, though the fact that applying a Legendre transform to a \gFM gives a \FM is an implicit link, via the structure of the Legendre field. Before considering the nature of the associativity equation in 2-dimensions we weaken this link further and consider more general structures.

\subsection{WDVV pairs}

By itself, one cannot form an associativity condition from a prepotential $F$, only a $(0,3)$-tensor from its third derivatives 
\[
c_{ijk} = \frac{\partial^3 F}{\partial t^i \partial t^j \partial t^k}\,.
\]
Without a symmetric $(2,0)$-tensor (the inverse metric) one cannot form a $(1,2)$-tensor and hence form a multiplication. This leads to:

\begin{defn}
A WDVV pair $\{F,\eta^\sharp\}$ consists of a prepotential $F$ and a symmetric $(2,0)$-tensor $\eta^\sharp$ such that the WDVV equations hold.
\end{defn}

\noindent To form the WDVV equations one only requires the existence of a symmetric $(2,0)$-tensor. In the above definition, the tensor $\eta^\sharp$ will be referred to as an inverse-metric even through, without a non-degeneracy assumption, its inverse - which would normally be the metric - may not exist. Note too that any conformal rescaling of $\eta^\sharp$ will also define a WDVV pair. Given $F$, the WDVV equation can then be seen as an overdetermined linear system for the components of the $\eta^\sharp$ in terms of the derivatives of $F\,.$ This will give the most general connection between these tensors, and this system is very easy to analyze in 2-dimensions.

But before that, we note that, given one WDVV pair, one may construct an infinite number of WDVV pairs.
\begin{prop}
    Given a WDVV pair $\{F,\eta^{ij}\}$ and an arbitrary non-zero vector field $\mathcal{E}$ then on defining
    \[
    g^{ij} = \eta^{ir}\eta^{js} c_{rst}\, {\mathcal{E}^t}
    \]
    then $\{F,g^{ij}\}$ is also a WDVV pair.
\end{prop}
\begin{proof}
By definition
\[
    F_{ijp}\, g^{pq} \, F_{qrs} = F_{ijp}\, \underbrace{(\eta^{pm} \, \eta^{qn}\, F_{mnt} \mathcal{E}^t)}_{g^{pq}} F_{qrs}
\]
and hence, on using the fact that the \lq seed\rq~pair 
$\{F,\eta^\sharp\}$ satisfies the WDVV equations, one can swap the $i$ and $r$ indices and this shows that 
$\{F,g^{ij}\}$ is also a WDVV pair.

\end{proof}

\noindent For example, if $\mathcal{E}$ is the Euler field then $g$ is the intersection form. From this one may form a recursion operator and generate a tower of structures. If one wants to guarantee some other structure is preserved, such as an $F$-manifold structure, then $\mathcal{E}$ must be an eventual identity \cite{davidstrDubduality}. This result is akin to the result in \cite{K1} which uses the matrix expansion $(1-A)^{-1}= 1+A+A^2 +\ldots\,.$

\subsection{Associativity in two dimensions}

We now restrict to two dimensions and analyse the structure of a $\{F\,,\eta^\sharp\}$-pair. Assuming this inverse metric has components (which may depend on the coordinates)
\[
\eta^{ij} = \left(\begin{array}{cc}
a & - b \\ -b & c \end{array}\right)
\]
then in two dimensions the WDVV-equations reduce to the single equation
\[
\left| \begin{array}{ccc}
F_{122} & F_{222} & a \\
F_{112} & F_{122} & b \\
F_{111} & F_{112} & c
\end{array}
\right|=0\,.
\]
One consequence of this is that it implies that the otherwise over-determined system
\[
\left(\begin{array}{cc}
F_{122} & F_{222}  \\
F_{112} & F_{122}  \\
F_{111} & F_{112}  
\end{array}
\right)
\left(\begin{array}{c} \alpha(t^1,t^2) \\ \beta(t^1,t^2) \end{array}\right)=
\left(\begin{array}{c} a \\ b \\ c \end{array}\right)
\]
will generically have a unique solution as the augmented matrix has rank 2 rather than 3 due to the associativity condition (for the lower-rank cases see \cite{AndreThesis}). So generically this determines a vector 
\[
e = \frac{1}{D} \left(\alpha(t^1,t^2) \frac{\partial~~}{\partial t^1} + \beta(t^1,t^2) \frac{\partial~~}{\partial t^2}\right)
\]
(where $D=ac-b^2$ is assummed to be non-zero) which turns out to be the identity for the multiplication. This result generalizes to higher dimensions \cite{AF}.

Restricting to the anti-diagonal case where $\eta=2 \, dt^1\, dt^2$ gives the associativity equation
\begin{equation}
F_{111} \, F_{222} - F_{112} \,F_{122} =0\, ,
\label{2dwdvv}
\end{equation}
which has been the central object of study throughout this paper. This Monge-Amp\`ere type-equation was studied in \cite{S96} where a Lax pair was constructed and the system reduced to a system of hydrodynamic type (with a single Riemann invariant, of Bianchi class VII). Writing (\ref{2dwdvv}) in the form
\[
\left|
\begin{array}{cc}
F_{11,1} & F_{22,1}\\
F_{11,2} & F_{{22,2}}
\end{array}
\right|=0
\]
implies there is a functional relationship
\[
F_{22}=f^{\prime\prime}(F_{11})
\]
for some arbitrary function $f\,.$ This reduces the order of the differential equation \cite{Kito} and also shows that any solution becomes, after the Legendre transformation generated by $\partial_{t^1}\,,$ the solution
\[
F=\frac{1}{2} {(t^1)^2 t^2 + f(t^2)}\,.
\]
Furthermore, if (\ref{2dwdvv}) holds then one may derive the identity
\[
\left(\begin{array}{cc}
F_{222} &  F_{112} \\ 
F_{122} & F_{111} 
\end{array}\right)
\left(\begin{array}{c}
\Delta_1 \\ \Delta_2
\end{array}\right) = \left(\begin{array}{c}
0 \\ 0
\end{array}\right)
\]
where $\Delta_1 = F_{111} F_{122} - F_{112}^2$ and $\Delta_2 = F_{222} F_{112} - F_{112}^2\,,$ and if these are non-zero then the inverse vector fields $\partial_1^{-1}$ and $\partial_2^{-1}$ exist. A simple analysis of these equation show - if the Frobenius algebra is to have no nilpotent elements - that both these inverses must exist.

\subsection{The solution space of the WDVV equations in 2-dimensions}

In this section we concentrate on the space of generalized Frobenius manifolds that reduce to the dihedral Frobenius manifolds. This analysis will also single out the special points in solution space identified above, but via a different mechanism.

Recall that in Example \ref{example51} the condition
\[
\widetilde{F}_{22} = \frac{1}{2} ( \widetilde{F}_{11} )^2
\]
immediately implied that $F_{22} = \frac12 (\tilde{t}^2)^2$ and hence after the application of the Legendre transformation one obtained the $I_2(3)$-Frobenius manifold. This is easily generalized. Consider the differential equation
\[
\widetilde{F}_{22} = \frac{1}{N-1} ( \widetilde{F}_{11} )^{N-1}
\]
which both implies the WDVV equation (\ref{2dwdvv}) is satisfied and that the Legendre-transformed solution is the dihedral $I_2(N)$ Frobenius manifold. Since
\[
E=\left(1-\frac{1}{N}\right) {\tilde{t}^1} \frac{\partial~~}{\partial {\tilde{t}}^1} + \frac{1}{N} {\tilde{t}^2} \frac{\partial~~}{\partial {\tilde{t}}^2}
\]
(and $\tilde{d}=1$), the scaling ansatz
\[
{\widetilde{F}} = (\tilde{t}^2)^{2N} g \left[\frac{\tilde{t}^1}{\tilde{t}^{N-1}} \right]
\]
reduces this to the second order nonlinear equation
\begin{equation}\label{reducedWDVV}
(N-1)^2 x^2 g^{\prime\prime}(x) - 3 N(N-1) x g^\prime(x) + 2N(2N-1) g(x) = \frac{1}{N-1} \left( g^{\prime\prime}(x) \right)^{N-1}\,.
\end{equation}
Thus any solution to this equation will generate a generalized Frobenius reducible, after a Legendre transformation to the $I_2(N)$-case. A leading order $g(x) = a_0 x^M + {\rm l.o.t.}$ analysis shows that no polynomial solution exists if $N\geq 5\,.$

\bigskip

\noindent{\underline{$N=2$}:}
This is a degenerate case as one obtains a cubic Frobenius manifold after the Legendre transformation. Equation (\ref{reducedWDVV}) becomes linear in this case and the full solution is
\[
g(x)= c_0 (1+x)^4 + c_1 (1-x)^4\,.
\]
\noindent This correspond to the generalized Frobenius manifold on the affine Coxeter group orbit spaces $\mathbb{C}^2/C_2$ (case 2) as found in \cite{orbit2}.

\noindent{\underline{$N=3$}:} From equating the highest power of $x$ in the differential to zero one obtain the equation $a_0(12 a_0 -1)=0$ and hence two cases. This lead to the two polynomial solution

\begin{eqnarray}
g(x) & = & \frac{27 a_1^4}{5} + 9 \,a_1^3 x + 9 a_1^2 \,x^2 + a_1\, x^3\,,\label{blue}\\
g(x) & = & \frac{a_2^2}{15} + a_2 x^2 + \frac{x^4}{12}\,.\label{green}
\end{eqnarray}

\noindent These correspond to the generalized Frobenius manifolds on the affine Coxeter group orbit spaces $\mathbb{C}/A_2$ and 
$\mathbb{C}/G_2$ respectively, as found in \cite{orbit1}.

\noindent{\underline{$N=4$}:}
From equating the highest power of $x$ in the differential to zero one obtain the equation $a_0(36 a_0^2 -1)=0$ and hence three cases. However the $a_0=0$ leads to a trivial solution. These lead to the two polynomial solution. 

\[
g(x) = \frac{a_1^3}{21} \pm \frac{a_1^2 x}{5} + a_1 x^2 \pm \frac{x^3}{6}\,.
\]

\noindent This correspond to the generalized Frobenius manifold on the affine Coxeter group orbit spaces $\mathbb{C}^2/C_2$ (case 1) as found in \cite{orbit2}.

\bigskip

One class of solutions - for all $N$ - can be found by considering the kernel of the Euler operator on the left-hand-side which annihilates the terms
\[
\displaystyle{x^\frac{2N}{N-1}\,,\qquad x^\frac{2N-1}{N-1}\,.}
\]
Deforming the second to $g(x) = c_1 x^\frac{2N-1}{N-1}  + c_2 x$ then provides a solution, since the power of $x$ on the right-hand-side $(\frac{2N-1}{2} - 2)(N-1)=1$ matches the remaining term on the left hand side (the first term in the solution is annihilated and the linear power in the second terms is preserved). In detail, this gives the solution
\begin{equation}
\widetilde{F} = c_1 (\tilde{t}^1)^\frac{2N-1}{N-1}\, \tilde{t}^2 + c_2 \tilde{t}^1 (\tilde{t}^2)^{N+1}\,.
\label{red}
\end{equation}
where the constants $c_1$ and $c_2$ are related by an inessential algebraic relation which depends on the normalizations used. On applying the Legendre transformation and inverting the change of variables show that this solution corresponds to the points $\alpha_1=0$ or $\alpha_2=0$. To see where the solution with $\alpha_1-\alpha_2=0$ and $\alpha_1+\alpha_2=0$ sit in the solution space of (\ref{reducedWDVV}) the phase-space portrait of the system will be analyzed. 

Here we present the $N=3$ case - the over-all structure remains the same for arbitrary $N\,.$ The scaling symmetry of the above equation allows one to introduce the scaling invariants
\[
f=\frac{g}{x^4}\,,\quad u = x \,\frac{g^\prime}{g}\,,\quad w = x^2 \,\frac{g^{\prime\prime}}{g}
\]
and on defining $v=x u^\prime$ and a change of independent variable $x=e^t$ reduces the equation to the dynamical system
\begin{eqnarray*}
    \dot{u} & = & v \,,\\
    \dot{v} & = & R(u,v)
\end{eqnarray*}
where 
\[
R(u,v)=\frac{21 u^4 - 2 u^5 + 2 v (45 + v) - 2 u^3 (39 + 5 v) + u^2 (119 + 65 v) - 
 u (60 + 139 v + 8 v^2)}{3 (10 - 7 u + u^2 + v)}\,.
\]
The critical points occur at the points $(0,0)\,,(1,0)\,,(5/2,0)\,,(3,0)$ and $(4,0)\,,$ theses values corresponding to the maximum and minimum powers of $x$ in the above explicit solutions.

\begin{figure}[ht]
    \centering
    \includegraphics[width=7cm]{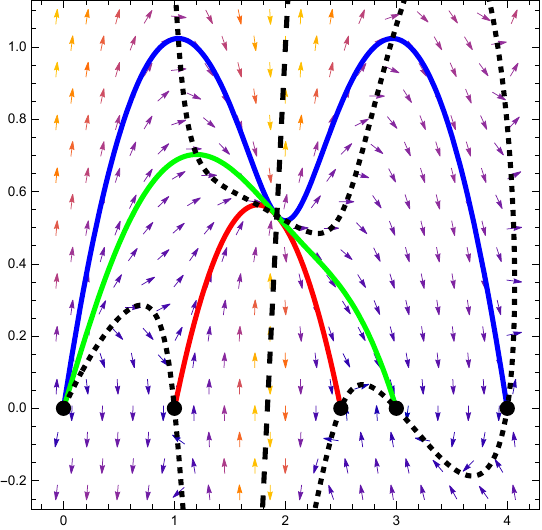}
    \caption{Phase space portrait for $N=3$ case, also showing the explicit, known, solutions}
    \label{fig:placeholder}
\end{figure}

Figure \ref{fig:placeholder} shows the phase space portrait for this system. Any solution to (\ref{reducedWDVV}) corresponds to a trajectory in this plane: the blue curve is the solution (\ref{blue}), the green the solution (\ref{green}) and the red the solution (\ref{red}). The curve consisting of small dots is where the numerator of $R(u,v)$ vanishes and the curve consisting of large dots is where the denominator vanishes. Note that all known solutions go through the singular point where these curves intersect. Furthermore, the known solutions leave/approach the critical points along eigenvector-directions. But recall that any solution corresponds to a point in the space $\mathbb{P}^1_{\tilde{d}=1}$ and the blue curve corresponds to the point $[1,1]$, the green curve to the point $[1,-1]$ and the red curve to the point $[1,0]$ or $[0,1]\,.$

This analysis suggests the following:

\begin{conj}

\phantom{text}

\noindent For arbitrary $N\in\mathbb{N}_{\geq 3}$ there exists trajectories:

\begin{itemize}

\item{} from $(0,0)$ to $(\frac{2N-2}{N-2},0)$ leaving/approaching these critical points along eigen-directions and corresponding to the point $[1,1] \in \mathbb{P}^1_{\tilde{d}=1}\,;$

\vskip 2mm

\item{} from $(0,0)$ to $(\frac{2N}{N-1},0)$ leaving/approaching these critical points along eigen-directions and corresponding to the point $[1,-1] \in \mathbb{P}^1_{\tilde{d}=1}\,.$

\vskip 2mm

\item{} from $(1,0)$ to $(\frac{2N-1}{N-1},0)$ leaving/approaching these critical points along eigen-directions and corresponding to the point $[1,0] {\rm~or~} [0,1]\in \mathbb{P}^1_{\tilde{d}=1}\,.$

\end{itemize}

\noindent In particular, these trajectories are heteroclinic orbits.

\end{conj}

\noindent The conjecture is true for $N=3\,,$ two-thirds true for $N=4$ (the blue and red curves are known explicitly) and one-third true for $N\geq 5$ (only the red curve is known explicitly). It should be possible to prove the existence of these heteroclinic orbits by a careful phase space analysis, but this would not show where they sit in the solution space.

\section{Inversion symmetry for generalized Frobenius manifolds}

So far we have concentrated on how Legendre transformation act on generalized Frobenius manifolds. Another symmetry of the WDVV equations were introduced in \cite{Dubrovin1996}: inversion symmetry. However the direct application of this to generalized Frobenius manifolds fails. To see this, and to set up the notation, we start with the theory of conformal transformations.

Since $e$ is a Legendre field, it has an associated potential $h_e$ and this may be used to define a conformal transformation, so
\begin{eqnarray*}
{}^I{\eta}(X,Y) & = & \frac{1}{h_e^2} \eta(X,Y)\,,\\
& = &  e^{2\phi}\eta(X,Y)\,.
\end{eqnarray*}
Then standard formulae from conformal geometry give
\begin{equation}
\begin{array}{ccc}
    {}^I{\nabla}_XY & = & \nabla_XY + (d\phi)(X)Y + (d\phi)(Y)X - \eta(X,Y) \nabla\phi\,,\\
    & = & \displaystyle{\nabla_XY - \frac{1}{h_e}\left\{
    \eta(e,X)Y + \eta(e,Y)X - \eta(X,Y)e\right\}}
\end{array}
\label{conformalconnection}
\end{equation}
and
\[
{}^I{Riem} = e^{2\phi} \, Riem + e^{2\phi} \eta \owedge H
\]
where $Riem$ is the $(0,4)$-Riemannian curvature tensor,
\[
H(X,Y) = \left( Hess\phi - d\phi \otimes d\phi + \frac{1}{2} |d\phi|^2\right)(X,Y)\,,
\]
and $\owedge$ is the Kulkarni–Nomizu product. A short calculation gives
\[
H(X,Y) = \left\{ -\frac{1}{h_e}  \eta(X\circ Y, \nabla_ee) + \frac{1}{2h_2^2} \eta(e,e) \eta(X,Y)\right\}\,.
\]
From this, the obstruction for $\widetilde{Riem}$ to be zero is
\[
\nabla_ee - \frac{1}{2 h_e} \eta(e,e)\, e =0\,.
\]
This immediate shows why the construction works for Frobenius manifolds but fails, in general, for generalized Frobenius manifolds. For Frobenius manifolds $\nabla_ee=0$ and, recalling that inversion symmetries are only defined in \cite{Dubrovin1996} if $\eta(e,e) = 0$, it follows that $\widetilde{\eta}$ is flat. While a Legendre transformation does not alter the spectrum of a Frobenius manifolds, inversion symmetries do. 

Using (\ref{conformalconnection}) it is easy to show that
\[
{}^I \mu(X) = \mu(X) + \left\{(d-1)X + \frac{1}{h_e} \eta(e,E)X\right\} + \frac{1}{h_e} \left\{ \eta(e,X)E - \eta(X,E)e \right\}
\]
and the middle term on the right vanishes by virtue of Lemma \ref{niceresult} (iv). Hence
\[
{}^I \mu(X) = \mu(X) +  \frac{1}{h_e} \left\{ \eta(e,X)E - \eta(X,E)e \right\}\,.
\]
In the standard coordinate for a Frobenius manifold this gives
\[
\tilde{\mu}_i = \begin{cases}
-1 + \mu_n &{\rm if~}i=1,\\
~~\mu_i &{\rm for~}i=2\,,\ldots\,,n-1,\\
+1-\mu_1 &{\rm if~}i=n.\\
\end{cases}
\]

For the above reasons, such an inversion symmetry cannot be applied to a generalized Frobenius manifolds - even if $\eta(e,e)=0,$ the unity field is not, in general, covariantly constant. To overcome this we adopt the following procedure. Starting with a generalized Frobenius manifold we apply a flat, homogeneous Legendre field to generate a Frobenius manifold. Under a mild algebraic condition on the Legendre field one can then apply the standard inversion symmetry to obtain a new Frobenius manifold. A new Legendre field can then be constructed from the original field to map to an inverted generalized Frobenius manifold. This is shown diagrammatically in the following:

\begin{equation*}
       \begin{tikzcd}[every cell/.append style={align=center}]
   \begin{tabular}{l}
   $\left\{\begin{array}{c}
{\rm generalized~Frobenius~~}\\
{\rm manifold:}
~ \{\eta\,,\nabla\}
\end{array}\right\}$
\\[-5pt]
        \arrow[dd, rightarrow, "\text{$\delta$ flat, homogeneous}"]
      \end{tabular} & $\overset{I(\delta)}{\joinrel\relbar\joinrel\relbar\joinrel\relbar\joinrel\relbar\joinrel\relbar\joinrel\rightarrow}$ & \begin{tabular}{c}
         $\left\{ 
         \begin{array}{c}
         {\rm generalized~Frobenius}\\
         {\rm manifold:~}
         ~ \{{}^{I(\delta)}\eta\,,{}^{I(\delta)}\nabla\}
         \end{array}
         \right\}$ 
         \\[-5pt]
          \arrow[dd, leftarrow, "\text{$\delta_I$ }"]
      \end{tabular}\\
        & &\\
  \begin{tabular}{c}
         $ \left\{
\begin{array}{c}
{\rm Frobenius~manifold:}\\
~\{\tilde\eta\,,\widetilde\nabla\}
\end{array}
\right\}$
        
      \end{tabular}  & $\overset{ {\rm~Inversion}}{\joinrel\relbar\relbar\joinrel\relbar\joinrel\relbar\joinrel\relbar\joinrel\relbar\joinrel\relbar\joinrel\rightarrow}$&\begin{tabular}{c} 
             $ \left\{
\begin{array}{c}
{\rm Frobenius~manifold:}\\
\{{}^I\tilde\eta\,,{}^I\widetilde\nabla\}
\end{array}
\right\}$ 
          
      \end{tabular} 
    \end{tikzcd}
\end{equation*}

\noindent Since $\delta^{-1}$ and $e$ are two Legendre fields defined in terms of the connection $\widetilde\nabla$ we have corresponding potentials $\tilde{h}_{\delta^{-1}}$ and $\tilde{h}_e\,.$

\begin{thm}
Let $\delta$ be a flat homogeneous Legendre field for a generalized Frobenius manifold with metric $\eta$ and connection $\nabla\,.$ 
For $\delta \neq e\,,$ the vector field
\[
{}^I \delta  = \tilde{h}_e \, \delta^{-1} - \tilde{h}_{\delta^{-1}} \, e
\]
is a Legendre field for ${}^I\widetilde{\nabla}\,.$ Moreover,
\[
{}^I \delta^\flat  = d \left(\frac{\tilde{h}_{\delta^{-1}}}{\tilde{h}_e}\right)\,.
\]
\end{thm}

\begin{proof} The proof is by direct computation. Recall that the potentials are defined by
\begin{eqnarray*}
    (d \tilde{h}_e)(X) & = & \widetilde{\eta}(e,X)\,,\\
    (d\tilde{h}_{\delta^{-1}})(X) & = & \widetilde\eta({\delta^{-1}},X)\,,
\end{eqnarray*}
and these are related to each other:
\[
\begin{array}{cclcl}
    (d\tilde{h}_e)(X) & = & \widetilde{\eta} (\delta\circ\delta^{-1},X)\,,&&\\
    & = & \widetilde{\eta}(\delta^{-1},X\circ \delta) & = &   (d\tilde{h}_{\delta^{-1}})(X\circ\delta)\,,\\
    & = & \eta(\delta, \delta\circ X) & = & (dh_\delta)(\delta\circ X)\,.\\
\end{array}
\]
To find a Legendre field ${}^I\delta$ one first needs the metric ${}^I\widetilde{\nabla}$ and connection ${}^I\widetilde{\nabla}\,.$ These are given by
\begin{eqnarray*}
{}^I\widetilde{\eta}(X,Y) & = & \frac{1}{({\tilde{h}_e})^2 } \,\tilde\eta(X,Y)\,\\
&&\\
{}^I\widetilde{\nabla}_XY & = & \widetilde\nabla_XY - \frac{1}{\tilde{h}_e} \left\{
    \tilde\eta(e,X)Y + \tilde\eta(e,Y)X - \tilde\eta(X,Y)e\right\}\,\\
\end{eqnarray*}
and the proof now is just the calculation of
\[
X \circ {}^I\widetilde{\nabla}_Y ({}^I\delta) -
Y \circ {}^I\widetilde{\nabla}_X ({}^I\delta)
\]
using these formulae. Since $\widetilde\nabla_X e = \delta^{-1}\circ \nabla_X\delta = 0$ (as $\delta$ is flat and homogeneous in order to map to a Frobenius manifold) the only obstruction for ${}^I\widetilde{Riem}$ to be flat is $\widetilde{\eta}(e,e) = 0$ or just $\eta(\delta,\delta)=0\,.$ Since $\delta$ just be a coordinate differential (say, $\delta = \frac{\partial~~}{\partial t^\sigma}$) this just requires $\eta_{\sigma\sigma}=0$. If the metric is anti-diagonal, this condition is automatically satisfied if the dimension is even, and excludes a single case if the dimension is odd. 

The final part is also a direct computation:
\begin{eqnarray*}
({}^I\delta)^\flat(X) & = & {}^I\tilde{\eta}({}^I\delta,X)\,,\\
& = & \frac{1}{ \tilde{h}_e } \widetilde{\eta}(\delta^{-1},X) - \frac{ \tilde{h}_{\delta^{-1}} }{(\tilde{h}_e)^2} \widetilde\eta(e,X)\,,\\
& = &  \frac{1}{ \tilde{h}_e } d \tilde{h}_{\delta^{-1}}(X) - \frac{\tilde{h}_{\delta^{-1}}  }{  (\tilde{h}_e)^2  } d\tilde{h}_e(X)\,,\\
& = & d \left(\frac{\tilde{h}_{\delta^{-1}}}{\tilde{h}_e}\right)(X)\,.
\end{eqnarray*}

\end{proof}

\noindent This theorem thus provides an inversion-like symmetry $I(\delta)$ between generalized Frobenius manifolds. Unlike the original inversion symmetry, it depends on a choice of a flat, homogeneous Legendre field (that is, a coordinate direction with $\eta(\delta,\delta)=0$) for the original generalized Frobenius manifold. The symmetry $I(\delta)$ does not change the charge; a short calculation shows that
\[
\mathcal{L}_E ({}^I\delta) = (1-d + \kappa) \, ({}^I\delta)
\]
from which ${}^Id=d$ follows on use of the results in Section \ref{LFsection}. It is interesting to note that while $I(\delta)$ is constructed using a suitable Legendre field, the charge and spectrum of the resulting generalized Frobenius manifold do not depend on properties of the field.

\section{Hurwitz spaces and the action of Legendre fields}

As mentioned in the introduction, one of the earliest examples of what is now called a (semi-simple) generalized Frobenius manifold had a simple Hurwitz space description where a `generalised' primary form -- often called \textit{horizontal differential} -- was used. Since applying a flat, homogeneous Legendre field to such a generalized Frobenius manifold results in a Frobenius manifold, it is a natural question to ask how Legendre transforms act at the level of Hurwitz space data.

We recall the main points of this construction (for full details see \cite{Dubrovin1996} and, for the specific context of generalized Frobenius manifolds, \cite{AlessandroThesis}). Given a Landau-Ginzburg superpotential $\lambda: \Sigma_g \rightarrow \mathbb{P}^1$ with some set of prescribed zeros and poles defining (a subset of) the Hurwitz space $H_{g;n}$, the metric and multiplication are given by the residue formulae
\begin{eqnarray}
 \label{eq:metricresidue}   \eta_\omega(X,Y) &=& \sum_{d\lambda=0} \Res \left\{
    \frac{X(\lambda) Y(\lambda)}{d\lambda} \, \omega^2\right\}\,,\\
    \eta_\omega(X\circ Y,Z) &=& \sum_{d\lambda=0} \Res  \left\{
    \frac{X(\lambda) Y(\lambda) Z(\lambda)}{d\lambda} \, \omega^2\right\}\,,
\end{eqnarray}
where $\omega$ is the primary differential. The pair $(\lambda,\omega)$ will be referred to as the Landau-Ginzburg model and the resulting object a Hurwitz Frobenius manifold. It is important to remark that such Frobenius manifolds are inherently semi-simple. 

An important object in the construction is the tautological bundle $C_{g;n}$ over $H_{g;n}\,$, whose fibre over $(\Sigma_g,\lambda)$ is given by the curve $\Sigma_g$ itself. This is the universal object for families of curves over $H_{g;n}\,$, and is therefore called \textit{universal curve} \cite{RomanoDoubleHurwitz}. Both the Hurwitz Frobenius manifold and its almost-dual extend to flat F-manifolds on the universal curve, endowing it with the structure of a bi-flat F-manifold \cites{dealmeida2025openhurwitzflatf,PS25}. The eventual identity $\EE$ coincides with the natural extension to the universal curve of the Euler vector field, as denoted by $D_E$ in \cite{RomanoDoubleHurwitz}. In particular, the key property that makes it the eventual identity between the two flat F-manifold structures is its action on the superpotential: $\Lie_\EE\lambda=\lambda\,$. See also \cite{AlessandroThesis} for further details.

Following \cite{RomanoDoubleHurwitz}, the structure induced by a horizontal, non-primary differential $\omega$ in general fails to fulfill the axiom regarding flatness of the unity, while homogeneity is preserved by the following class:
\begin{defn}
    A horizontal differential $\omega$ is \textit{homogeneous} if it is an eigenvector of $\Lie_\EE\,$:
    \[
    \Lie_\EE\omega=d_\omega\,\omega\,.
    \]
    The eigenvalue $d_\omega$ is called the \textit{degree} of $\omega\,$.
\end{defn}
The Euler vector field $E$ is a homothetic vector field\footnote{That is, it is a conformal Killing vector field with constant conformal factor.} for the flat  metric $\eta_\omega$ associated to a homogeneous horizontal differential $\omega\,$. As shown e.g. in \cite{Rej23}, in fact, in such case it is easy to see that:
\[
\Lie_E\eta_\omega=(1+2d_\omega)\,\eta_\omega\,.
\]
As a consequence, the corresponding (generalised) Frobenius manifold has charge $d=1-2d_\omega\,$.

\begin{prop} Let $(\lambda,\omega)$ be the Landau-Ginzburg model for a Hurwitz Frobenius
manifold, and $\delta$ be a Legendre field on the underlying Hurwitz space $H_{g;n}$. A Landau-Ginzburg model for the Legendre-transformed Frobenius manifold is $(\lambda,\widetilde{\omega})$, with:
\begin{equation}
\widetilde\omega := (\Lie_\delta \lambda)\, \omega\,.
\label{LGtrans}
\end{equation}
\end{prop}

\begin{proof}
From \cite{SS}, the tensor defining the multiplication does not change under a Legendre transformation, only the metric, and hence the flat coordinates for this metric. Recall that
\begin{eqnarray*}
    \widetilde\eta(X,Y) & = & \eta(X\circ \delta,Y\circ \delta) \,, \\
    & = & \sum_{d\lambda=0} \Res  \left\{
    \frac{(\delta \circ X)(\lambda) (\delta\circ Y)(\lambda) }{d\lambda} \, \omega^2\right\}\
\end{eqnarray*}
Hurwitz Frobenius manifolds are always semi-simple with canonical coordinates being defined by $u^i=\lambda(\alpha^i)$ where $\lambda^\prime(\alpha^i)=0\,.$ Using these coordinates
\begin{eqnarray*}
    \widetilde\eta(X,Y) & = & \sum_{d\lambda=0} \Res 
    \left\{ \frac{\delta^i \delta^j \partial_i \lambda \partial_j\lambda}{d\lambda} \, \omega^2\right\}\, X^i\, Y^j\,,\\
    & = & \sum_{d\lambda=0} \Res 
    \left\{ \frac{ \partial_i \lambda \partial_j\lambda}{d\lambda} \, (\delta(\lambda)\,\omega)^2\right\}\, X^i\, Y^j\,,\\
& = & \sum_{d\lambda=0} \Res  \left\{
    \frac{X(\lambda) Y(\lambda)}{d\lambda} \, \widetilde{\omega}^2\right\}\,.
\end{eqnarray*}
This uses the expansion around the critical point $p=\alpha^i\,,$
\[
\lambda(p) = u^i + \tfrac{1}{2} (p-\alpha^i)^2 \lambda^{\prime\prime}(\alpha^i) + O(p-\alpha^i)^3 
\]
so
\[
\left.\delta(\lambda)\right|_{p=\alpha^i} = \delta^i\,.
\]
This shows that the Legendre transformed structures come from the same Landau-Ginzburg superpotential as the original structure but with a transformed primary differential given by (\ref{LGtrans})\,.
\end{proof}

\begin{oss}
    Equation (\ref{LGtrans}) has a curious history. It appears, for the original Legendre transformations (so for flat-homogeneous Legendre fields) in \cite{Dubrovin1996} as the unnumbered equation between equation (5.82) and (5.84) - there is no equation with the given number (5.83)!
\end{oss}
Note that, in the proof, there is no requirement that the Legendre field should be flat or homogeneous, and so this transformation of the differential $\omega\mapsto\widetilde{\omega}$ provides a symmetry of the WDVV-equations. Having established the way in which the Landau-Ginzburg model transforms under Legendre transformation, one should then check that homogeneous Legendre fields preserve the class of generalised Frobenius manifolds, and that the charge changes according to \cref{prop:spectrumhomogeneousLF}.
\begin{cor}
    Let $(\lambda,\omega)$ be a Landau-Ginzburg model, $\omega$ being homogeneous of degree $d_\omega\,$. If $\delta$ is a homogeneous Legendre field with degree $\kappa\,$, then the Legendre-transformed horizontal differential $\widetilde{\omega}$ is homogeneous with degree:
    \[
    d_{\widetilde{\omega}}=d_\omega+\kappa+1\,.
    \]
    As a consequence, the charge of the Legendre transformed structure is $\widetilde{d}=1-2d_{\widetilde{\omega}}=d-2(\kappa+1)\,$, in accordance with \cref{prop:spectrumhomogeneousLF}.
\end{cor}
\begin{proof}
    This is the result of a direct calculation:
     \[
    \begin{aligned}
        \Lie_\EE\widetilde{\omega}&=\Lie_\EE(\Lie_\delta\lambda)\,\omega+(\Lie_\delta\lambda)\,\Lie_\EE\omega\\
        &=(\Lie_\delta\lambda)\,\omega+\bigl(\Lie_{\comm{E}{\delta}}\lambda\bigr)\omega+d_\omega\,\widetilde{\omega}\\
        &=\bigl(1+\kappa+d_\omega\bigr)\,\widetilde{\omega}\,.
    \end{aligned}
    \]
    The vertical component of $\comm{\EE}{\delta}$ does not contribute because, under the assumption that $\omega$ is homogeneous, it produces a term proportional to $\dd\lambda\,$, which does not change the metric determined by $\widetilde{\omega}$ from \cref{eq:metricresidue}.
\end{proof}

 More details of this construction, including the action of Legendre transformation on the bi-differentials that live on Hurwitz spaces, may be found in \cite{AlessandroThesis}. Closely related results have recently appeared in \cite{Ma26}. We end with some examples.

\begin{ex}
    The generalized Frobenius manifold (\ref{nonFMsolution}) from the introduction was shown in \cite{BCR} to have the Hurwitz space description in terms of the pair
\begin{eqnarray*}
    \lambda(p) & = & \frac{p(p-e^v)}{p-e^{v+w}}\,,\\
    \omega & = &\tfrac{1}{p}\dd{p}\,,
\end{eqnarray*}
with flat coordinates $t^1 = e^v(e^w-1)\,, t^2 = v+w\,.$ Applying the Legendre transformations generated by $\delta = \partial_{t^1}$ gives
\begin{eqnarray*}
    \widetilde{\omega} & = & (\partial_{t^1}\lambda) \, \omega \,,\\
     & = & \frac{p}{p-e^{t^2}} \, \frac{1}{p}\dd{p}\,
\end{eqnarray*}
and hence $\tilde{p} = \log(p-e^{v+w})\,.$ In terms of this variable the new pair becomes
\begin{eqnarray*}
    \lambda(\tilde{p}) & = & e^{\tilde{p}} + \tilde{s}^1 +  \tilde{s}^2 e^{-\tilde{p}}\,,\\
    \widetilde{\omega} & = & \dd \tilde{p}
\end{eqnarray*}
which is the standard Hurwitz space data for the Frobenius manifold $FM_4$ on the extended affine Weyl orbit space $\mathbb{C}^2/A_1^{(1)}\,.$
\end{ex}

\begin{ex}
    The following is a two-dimensional version of the example appearing in \cite{RomanoThesis}. Consider the pair
\begin{eqnarray*}
    \lambda(p) & = & p\left(p^2 + t^2 p +t^1 + \tfrac{1}{6} t^2  \right)\,,\\
    \omega & = & \frac{1}{p}\dd{p}\,.
\end{eqnarray*}
Applying the Legendre transformation generated by $\delta = \partial_{t^1}$ gives
\begin{eqnarray*}
    \widetilde{\omega} & = & (\partial_{t^1}\lambda) \, \omega \,,\\
     & = & p \, \frac{1}{p}\dd{p}\,
\end{eqnarray*}
and hence $\tilde{p} = p\,.$ In terms of this variable the new pair becomes
\begin{eqnarray*}
    \lambda(\tilde{p}) & = & {\tilde{p}^3} + \tilde{s}^2 {\tilde{p}}+  \tilde{s}^1 \,,\\
    \widetilde{\omega} & = & \dd \tilde{p}
\end{eqnarray*}
which is the standard Hurwitz space data for the Frobenius manifold $FM_1$ on the Coxeter orbit space $\mathbb{C}^2/I_2(3)\,.$ In the last step a shift $\tilde{p} \rightarrow {\tilde{p}} + \mathrm{constant}$ has been used to put the superpotential into a standard form.

This all generalizes \cite{AlessandroThesis}, with the homogeneity condition giving a triangular change of variables between the coefficients of powers of $p$ in the superpotential and the flat coordinates. Thus 
\begin{eqnarray*}
    \lambda(p) & = & p\left(p^n + t^n p^{n-1} + \dots +t^1 + f(t^2\,,\ldots\,,t^{n}) \right)\,,\\
    \omega & = & \tfrac{1}{p}\dd{p}\,.
\end{eqnarray*}
and so taking the Legendre field $\delta = \partial_{t^1}$ maps this superpotential and primary differential to the standard form of the Frobenius $A_n$-superpotential
\begin{eqnarray*}
    \lambda(\tilde{p}) & = & {\tilde{p}^{n+1}} + \tilde{s}^n {\tilde{p}^{n-1}}+  \ldots +\tilde{s}^1 \,,\\
    \widetilde{\omega} & = & \dd \tilde{p}\,.
\end{eqnarray*}
This provides a map between the $(A_n,\omega_n)$-superpotential/differential developed in \cite{orbit2} and the classical superpotential/differential coming from the unfolding of the $A_n$-singularity.
\end{ex}

\bigskip

Clearly there is much scope to develop this in more generality starting with superpotentials and primary forms of the form (also called $(M,N)$-Ablowitz-Ladik structures)
\begin{eqnarray*}
\lambda(p) & = & p \frac{\prod_{i=1}^M (p-z_i)}{\prod_{j=1}^N (p-p_j)}\,,\\
\omega & = & \frac{1}{p}\dd{p}
\end{eqnarray*}
and further extension with multiplicities can also be developed.

\section{Comments}

In \cites{FKS,LeoThesis} flat Legendre fields were used to construct a map from certain rational $\bigvee$-systems to trigonometric $\bigvee$-systems, and this result can now be seen as an example of Theorem 3.1. Starting with the Frobenius manifold structure on the orbit space $\mathbb{C}^n/W$ for $W=A_n\,,B_n\,,D_n$ the almost dual Frobenius manifolds are basic examples of $\bigvee$-systems (where the $\bigvee$-system is just the root systems $R_W)$. Choosing an appropriate flat vector field $\delta^\star$ and applying the corresponding Legendre transformation resulted in certain trigonometric $\bigvee$-systems, and one question that arose from this work was whether there was a corresponding Frobenius-manifold type object whose dual object was given in terms of this trigonometric $\bigvee$-system. Theorem 3.1 identifies this object as a generalized Frobenius manifold:

\bigskip

\begin{equation*}
       \begin{tikzcd}[every cell/.append style={align=center}]
   \begin{tabular}{c}
   $\left\{\begin{array}{c}
{\rm Coxeter~}\\
{\rm orbit~space}
\end{array}\right\}$
\\[-5pt]
        \arrow[dd, leftrightarrow, "\text{duality}"]
      \end{tabular} & $\overset{\delta^{-1} {\rm ~flat}}{\leftarrow\joinrel\relbar\joinrel\relbar\joinrel\relbar\joinrel\relbar\joinrel\relbar\joinrel\relbar\joinrel\relbar\joinrel}$ & \begin{tabular}{c}
         $\left\{ 
         \begin{array}{c}
         {\rm Generalized~Frobenius}\\
         {\rm manifold~}
         (\widetilde{gFM})
         \end{array}
         \right\}$ 
         \\[-5pt]
          \arrow[dd, leftrightarrow, "\text{generalized duality}"]
      \end{tabular}\\
        & &\\
  \begin{tabular}{c}
         $ \left\{
\begin{array}{c}
{\rm Rational~}\\
{\rm \bigvee-{\rm system}}
\end{array}
\right\}$
        
      \end{tabular}  & $\overset{\delta^\star {\rm~flat}}{\relbar\joinrel\relbar\joinrel\relbar\joinrel\relbar\joinrel\relbar\joinrel\relbar\joinrel\relbar\joinrel\rightarrow}$&\begin{tabular}{c} 
             $ \left\{
\begin{array}{c}
{\rm Trigonometric~}\\
{\rm \bigvee-{\rm system}}
\end{array}
\right\}$ \\
          
      \end{tabular} 
    \end{tikzcd}
\end{equation*}

\bigskip

\noindent While $\widetilde{gFM}$ certainly exists, how its geometry is related to that $\mathbb{C}^n/W$ requires more work. Calculations suggest that $\widetilde{gFM}$ should be the affine orbit spaces constructed recently in \cites{orbit1,orbit2}. But this construction also raises other questions. In \cite{FKS} the coordinate-inversion problem could be done explicitly only for the Coxeter groups $W=A_n\,,B_n\,,D_n$ and only for certain special directions which should correspond to specific choices of the constants $\alpha_i\,.$ How these are related to the Coxeter geometry is unclear. Furthermore, this rational-to-trigonometric map was used for non-Coxeter $\bigvee$ systems for which there is no corresponding Frobenius manifold. How these fit into this picture is unknown, but certain $\bigvee$-systems arise on restriction of a Coxeter configuration onto discriminant submanifolds \cite{FV}. In \cite{SS} it was shown that the parallel component of a Legendre field when restricted to a discriminant submanifold was a Legendre field for the induced connection/metric on the submanifold. Thus it may be possible to extend this theory to cover maps between certain $\bigvee$-systems to certain trigonometric $\bigvee$-systems.

A final question is how the monodromy group and orbit space structures are related. In examples one gets a map between Frobenius manifolds on Coxeter group orbit spaces, and generalized Frobenius manifolds on affine-Weyl orbits spaces, and it would be interesting to see how this theory could be developed to provide a direct map between the corresponding monodromy groups and hence a map between these orbit spaces.

\noindent{\bf Acknowledgements}
 A. B. was supported by a Maclaurin Scholarship provided by the School of Mathematics and Statistics, University of Glasgow and L. K. by a Scholarship provided by the College of Science and Engineering, University of Glasgow. A. P. was supported by Ph.D. studentships from the EPSRC Doctoral Training Partnership (EP/W524359/1). We would also like to thank Paolo Lorenzoni for various useful conversations.

\medskip

\noindent{\bf Conflict of Interest:} On behalf of all authors, the corresponding author states that there is no conflict of interest.

\medskip

\noindent{\bf Data Availability:} Not applicable as article has no associated data.

\medskip

\noindent{\bf Use of software:} AI was used to proof-check the paper. Mathematica was also used to produce Figure 5.

\phantomsection
\addcontentsline{toc}{section}{References}
\bibliography{biblio} 

\end{document}